\documentclass[11pt]{article}
\usepackage{amssymb,epsf,latexsym, amsmath}
\makeatletter
\@addtoreset{figure}{section}
\def\thefigure{\thesection.\@arabic\c@figure}
\def\fps@figure{h, t}
\@addtoreset{table}{bsection}
\def\thetable{\thesection.\@arabic\c@table}
\def\fps@table{h, t}
\@addtoreset{equation}{section}
\def\theequation{\thesection.\arabic{equation}}
\makeatother

\def \eps{{\varepsilon}}
\def \n{{\eta}}
\def \kk{{\kappa}}
\def \kappa{{\nu}}
\def\la{\langle}
\def\ra{\rangle}
\def \n{\eta}
\def \tagm{\tilde{A}^\gamma}
\def \ta{\tilde{A}}
\def \b{{\mathcal B}}
\def \d{{\mathcal D}}
\def \t{{\mathcal T}}
\def \l{{\mathcal L}}
\def \lki{\Lambda_{ki}}
\def \m{{\mathcal M}}
\def \o{{\mathcal O}}
\def \r{{\mathbb R}}
\def \h{{\mathbb H}}
\def \q{{\mathbb Q}}
\def \p{{\mathbb P}}

\def \n{{\mathbb N}}
\def \v{{\mathbb V}}
\def \gd{{\Gamma^\delta}}
\def \G{G_{ik}^\kappa}
\def \Gx{{G_{ik}^\kappa}^1}
\def \Gy{{G_{ik}^\kappa}^2}
\def \s{{\mathcal S}^\eps}
\def \sd{{\mathcal S}^\eps_\delta}
\def \si{{S^{-1}}}
\def \N{{ \Gamma (N^\eps|_{\cup_{i=1}^s U^{J-3}_i}) }}

\def \tb{\overline{\tau} }
\def \sb{\overline{S} }
\def \sbk{\overline{S^\kappa} }
\def \clm{\overline{M} }
\def \sk{{ S^\kappa  }}
\def \uk{ u^\kappa  }
\def \lm{{ \rightarrow }}

\def \tbn{\overline{T}^n  }
\def \tbnk{\overline{T}^n_k  }
\def \tbni{\overline{T}^n_i  }
\def \tbbn{\overline{\overline{T}}^n  }
\def \tbbnk{\overline{\overline{T}}^n_k  }
\def \tbbni{\overline{\overline{T}}^n_i  }
\def \wu{W^{\bar{u}} }
\def \wuh{W_h^\rho }
\def \wup{W^{\bar{u},\rho} }
\def \nup{N^{\bar{u},\rho} }
\def \wuph{W_h^{\bar{u},\rho} }
\def \Puph{\Phi_h^{\bar{u},\rho} }

\allowdisplaybreaks

\begin{document}

\title {Persistence of Invariant Manifolds for Nonlinear PDEs}
\author{
Don A. Jones\thanks{Research partially supported by
CHAMMP.  email: dajones@schwinn.la.asu.edu}\\
Department of Mathematics\\
Arizona State University\\
Tempe, Arizona 85287-1804 USA
\and Steve Shkoller
\thanks{Research partially supported by the Cecil and Ida M. Green Foundation
and DOE. email: shkoller@cds.caltech.edu}
\\ Center for Nonlinear Studies
\\ Los Alamos National Laboratory  MS-B258\\
 Los Alamos, NM 87545 USA
}
\date{September 1996; this version, August 7, 1997\\
To appear in Studies in Applied Mathematics}

\maketitle

\newcommand{\Dt}{\Delta t}
\newcommand{\mm}{\mbox{$\mathcal M$}}
\newcommand{\ttt}{\mbox{$\mathcal T$}}
\newcommand{\hhh}{{\mathcal H}}
\newcommand{\aaa}{\mbox{$\mathcal A$}}
\newcommand{\dd}{\mbox{$\mathcal D$}}
\newcommand{\oo}{\mbox{$\mathcal O$}}
\newcommand{\rr}{\mbox{$\mathcal R$}}
\newcommand{\vv}{\mbox{$\mathcal V$}}
\newcommand{\bu}{{\bar u}}
\newcommand{\by}{\bar Y}
\newcommand{\rrr}{I\!\!R}
\newcommand{\dt}{\Delta t}
\newcommand{\R}[1]{$(\ref{#1})$}

\renewcommand{\theequation}{\thesection.\arabic{equation}}
\newtheorem{thm}{Theorem}[section]
\newtheorem{lemma}[thm]{Lemma}
\newtheorem{prop}[thm]{Proposition}
\newtheorem{cor}[thm]{Corollary}
\newtheorem{cond}[thm]{Condition}
\newtheorem{rem}[thm]{Remark}
\newtheorem{defn}[thm]{Definition}

\nopagebreak
\begin{abstract}
We prove that under certain stability and smoothing properties of
the semi-groups generated by the partial differential equations that
we consider, manifolds left invariant by these flows persist under
$C^1$ perturbation.  In particular, we extend well known 
finite-dimensional results to the setting of an infinite-dimensional
Hilbert manifold with a semi-group that leaves a submanifold invariant.
We then study the persistence of global unstable manifolds of
hyperbolic fixed-points, and as an application consider the 
two-dimensional Navier-Stokes equation under a fully discrete
approximation.  Finally, we apply our theory to the persistence
of inertial manifolds for those PDEs which possess them.
\end{abstract}

\tableofcontents

\section{Introduction}
We consider a nonlinear partial differential equation which takes the
form of an evolution equation having infinite-dimensional Riemannian
manifold $X$ as its configuration space. We assume that this PDE generates a $C^1$ 
nonlinear semi-group $S(t)$ which leaves a $C^r$ differentiable compact 
submanifold with 
boundary $\clm$ either overflowing or inflowing invariant, and examine
the persistence of this submanifold under small $C^1$ perturbations
of $S(t)$ in a neighborhood of $\clm$. Specifically, we prove that if
a semi-group $\sk(t)$ is within a $\kappa$-diameter $C^1$ ball about
$S(t)$, and if
the submanifold $\clm$ is normally hyperbolic, attracting trajectories
at an exponential rate, with a rate of attraction in the normal
direction greater than in the tangential direction, a 
continuously differentiable invariant 
submanifold with boundary of the same type $\overline{M_\kappa}$ exists
for $\sk(t)$ and converges to $\clm$ as $\kappa$ tends to zero.
Our result generalizes the persistence theory for manifolds
invariant to finite-dimensional vector fields that has been studied
by Moser, Sacker, Hirsh, Pugh, \& Shub, and Fenichel (see \cite{Moser}, 
\cite{Sac}, \cite{HPS} and \cite{Fen}).  Unlike the flows generated by 
differentiable finite-dimensional vector fields, however, the 
semi-group $S(t)$ does not, in
general, possess a bounded inverse and hence, for fixed $t$, the map
$S(t)$ is not a diffeomorphism.  In particular, its inverse is only
defined on the invariant manifold.  Furthermore, the configuration space
is not locally compact, so that one must strongly rely on the $C^1$ 
regularity of the map together with the compactness of $\clm$ to establish
important local bounds.

Recently, Bates, Lu, \& Zeng (see \cite{BLZ}) have generalized the
finite-dimensional persistence theory of Hirsh, Pugh, \& Shub to 
semiflows in Banach spaces; in our notation, their configuration space 
$X$ is a Banach vector space.  They consider compact manifolds $M$ 
without boundary as having center-type flow, and are able to prove the
persistence result for both stable-center and the more difficult
unstable-center case, the latter being quite complicated due to the
general lack of backwards-in-time flow.  Herein, our main interest
is in applying the general persistence theory to global unstable 
manifolds and inertial manifolds associated with nonlinear partial
differential equations, so we will restrict our attention to the case of
stable flow in the infinite-dimensional Riemannian manifold.

We remark that our choice of a Hilbert manifold for the configuration space 
may be attributed to the following reasons.  First,
we must require the existence of a second order vector field or spray on $X$
which generates the geodesic flow.  
In general, Banach spaces are not separable so 
that an infinite-dimensional manifold modeled on a Banach space may
not be paracompact and hence may not have global 
partitions of unity.  Since it is essential for us to be able to 
construct tubular neighborhoods 
about our embedded submanifold $\clm$ (we describe these in the next section),
we must require geodesics on $X$ to exist.  In the Hilbert manifold 
setting, the existence of such geodesic flows is equivalent to the
existence of a Lagrangian vector field associated with the kinetic
energy Lagrangian.  The Riemannian structure generates a natural 
symplectic form on $T^*X$ which is only weakly nondegenerate, meaning that
the associated linear map taking vector fields to $1$-forms is generally
not surjective.  It is, nevertheless, possible to confirm that the Hilbert
manifolds common to most applications do indeed possess geodesics, and herein
we shall restrict our attention to such examples. A classical example
of a configuration manifold which does not have linear structure arises in fluids
applications.
If one is interested in perfect incompressible
fluid flow governed by the Euler equations, then the approriate  
configuration space is the closure
of the volume preserving diffeomorphisms under a Hilbert space norm, a 
Hilbert manifold (see \cite{EM}).  
Our examination of the persistence problem for this PDE shall appear in a
future work.  Another interesting example is the Sine-Gordon equation
whose fields take values in ${\mathbb S}^1$.  One could argue that the
analysis of this PDE could be trivially done by working in the universal
covering space, however, if one is interested in the perturbation of
invariant manifolds, and wishes to work in the linear space supplied by
the cover, one would then have to prove that the perturbed structure
remains an invariant manifold after projecting from the covering space back
to the quotient.

Another reason for choosing a Hilbert space structure
is that we may reduce the normal bundle of the invariant manifold $\clm$
to a Hilbert or unitary bundle, wherein the transition maps are unitary.
This is essential for the local analysis that we perform.

Recall that an overflowing (inflowing) invariant manifold with 
boundary $\clm$ satisfies $S(t)u \in \clm$ for any $u\in\clm$ for all 
$t<0$ ($t>0$).  This means that the infinite-dimensional vector
field defined by the partial differential operator points outward
(inward) on the boundary of the manifold.
The existence of such an  overflowing (inflowing) invariant manifold 
$\overline{M_\kappa}$ for the perturbed  semi-group $\sk (t)$ is 
established by standard contraction mapping arguments. Specifically, 
after diffeomorphically identifying the infinite-dimensional 
manifold $X$ with the normal bundle over $\clm$ in a neighborhood 
$V$ of $\clm$, we search for $\overline{M_\kappa}$  in the space of
sections of this normal bundle.  The invariant manifold is the image of
the particular section which is the fixed-point for each time-$t$ map 
$\sk (t)$.  We remark that Fenichel's proof in finite dimensions  in
its local form essentially works for the infinite-dimensional vector
bundle setting with minor modifications, especially for showing that
the persisting manifold is Lipschitz.  We provide a different proof,
however, that the persisting manifold  is actually continuously
differentiable.

As an application of our persistence theorem, we show that the
unstable manifold of a hyperbolic fixed-point of a 
PDE satisfies the conditions sufficient for persistence so that the
approximating semi-group possesses a continuously differentiable
unstable manifold of its hyperbolic fixed-point which uniformly
converges to the unstable manifold of the unperturbed system
(see \cite{Mar} for a proof using the deformation method for
perturbation of the linearized system).
We may then deduce, following the work of \cite{Hump}, that the global
unstable manifold of the  obtained stationary solution, defined by
evolving the perturbed local unstable manifold forward in time, is
lower semi-continuous.   In case the global
attractor is the closure of the unstable manifolds (of overflowing
manifolds) such as with gradient systems, we
find that the  attractor is lower semi-continuous as well.

As the majority of nonlinear PDEs generate semi-groups which have no
explicit representation,  it is of great interest to consider
numerical schemes as the $C^1$ close approximations of the semi-group.
In fact, our motivation for perturbing the semi-group rather than the
PDE may primarily be attributed to the inability of numerical schemes
to approximate infinite-dimensional vector fields in a $C^1$ sense.
As an application, we show that the  unstable manifold of the 
hyperbolic stationary solution to the 2D Navier-Stokes equation
persists under  a fully discrete numerical approximation.

Next we apply our theory to the persistence of inertial manifolds
for those PDEs that possess them. An inertial manifold for a PDE is a
{\em smooth} finite-dimensional, exponentially attracting, positively
invariant manifold containing the global attractor (see \cite{FST1},
\cite{FST}). Subsequent works extended the existence results to more
general equations and provided alternative methods of proof (see for
example \cite{RT} and the references therein). Most existence proofs
require a gap condition to hold in the spectrum of the linear
term. As shown in \cite{RT}, this gap condition also insures that the
inertial manifold is normally hyperbolic.
Thus, our persistence theory applies to PDEs satisfying this gap
condition; however, since our theory is local, we may only conclude
that the perturbed system possesses an inflowing manifold that is
close to the inertial manifold. In particular, the inflowing manifold
for the perturbed system may not globally exponentially attract all
trajectories, as often occurs when discretizing dissipative systems.
We include cases where the inflowing manifold is an
inertial manifold for the perturbed system. 

Our applications complement existing results relating the long-time
dynamics of PDEs to the long-time dynamics of their numerical
approximations as the time and space mesh are refined. In
\cite{CFT}, \cite{Ti1} sufficient conditions are found on a
Galerkin scheme approximating the Navier-Stokes equations (NSE) to infer the
existence of a nearby stable stationary solution for the NSE from the
apparent stability of the time-dependent Galerkin approximate
solutions. Stable periodic orbits are studied in \cite{Ti3}, while
hyperbolic periodic orbits can be found in \cite{AD2}. Other aspects of the large-time
behavior of the Galerkin approximation to the exact solutions of the
NSE are considered in  \cite{HR}; subsequent results examined
the behavior of solutions near equilibria of PDEs under numerical
perturbation (see \cite{Stuart}, \cite{AD}, \cite{LSS}). 
For more general global attractors, see \cite{HLR} for sufficient
conditions  on upper semi-continuity.
For persistence of inertial manifolds under numerical approximation,
see \cite{FST}, \cite{FSTi}, \cite{DG}, \cite{FTi}, \cite{JS1}, and
\cite{JSTi}.

We structure the paper in the following way.  In Section 2,
we prove that if a nonlinear PDE has an overflowing invariant manifold
that is normally hyperbolic and whose semi-group satisfies certain
regularity conditions, then a continuously differentiable overflowing
invariant manifold persists under $C^1$ perturbation of the semi-group.
In Section 3, we state the general form  of the nonlinear partial
differential equation and define the domain of the semi-group which it
generates.  In Section 4, we prove that the global unstable manifold
of a hyperbolic fixed-point (of the PDE described in Section 3) 
persists as a continuously differentiable overflowing invariant
manifold. We apply this result to the examination of the two-dimensional
Navier-Stokes equation in Section 5, and finally, in Section 6,
we prove that inertial manifolds persist as well.

\section{Construction of the perturbed manifold}

In this section, we generalize Fenichel's proof of the persistence
of overflowing invariant manifolds to the setting of an infinite
dimensional Riemannian manifold.  Much of the proof goes through
essentially unchanged, particularly the portion which shows that the persisting
manifold is Lipschitz continuous; however, we provide a complete
exposition in order to account for the changes required in infinite
dimensions, as well as to clarify and generalize some of the details
in the technical lemmas.  As for the final result that the persisting
manifold is actually $C^1$, we take a different approach than Fenichel,
and show that the persisting manifold is a limit of a $C^1$ Cauchy
sequence.
\subsection{Infinite-dimensional geometry and sufficient conditions
for persistence}
\label{theory}

Let $(X,g)$ be an infinite-dimensional $C^r$
Riemannian manifold modeled on a
Hilbert space $(\h,\la \cdot, \cdot \ra )$, where $r\ge 3$. 
We consider a $C^1$ injective map $S:X\lm X$ which leaves
a submanifold with boundary  $\clm$ negatively 
invariant and has an inverse defined on $\clm$.
In particular, we assume the existence of a 
compact, connected, and embedded  $C^r$
submanifold $\clm = M \cup \partial M$ such that $S^{-1}(\clm) 
\subset \clm$.  We further assume that $\clm \subset M'$, where 
$M'$ is also a $C^r$ connected and embedded submanifold of $X$
having compact closure, which is the same dimension as $M$ and
satisfies $\si (M') \subset \clm$.  This gives each $x \in \clm$
an $M'$-open neighborhood, and thus allows us to avoid the complications
associated with charts defined on domains with boundary.
We assume that the differential of $S^{-1}|_{M'}$ is only defined in 
the tangent
bundle of $M'$ rather than in the entire tangent bundle of $X$ over
$M'$.

We consider
the tangent bundle $(T(X),\pi)$ restricted to $\pi^{-1}(M')$ and use
its Riemannian structure to induce the splitting
\begin{equation}
T(X)|_{M'} = T(M') \oplus N, \nonumber
\end{equation}
where $N$ is the $C^{r-1}$ normal bundle (any transverse bundle).  It 
is implicit in our compactness assumption that $\overline{M'}$ is 
finite dimensional, and for
concreteness, we fix the dimension to be $m$.  Due to the local 
trivialization of the
vector bundles, there exists an open covering  consisting of 
$M'$-open sets ${\mathcal U}=\{U_i, i\in {\mathcal I} \}$ of $\clm$,
such that for each $U_i \in {\mathcal U}$, we have the vector bundle 
morphisms
\begin{equation}
\tau^{M'}_i: T(M')|_{U_i} \lm U_i \times \p, \ \ \
\tau^N_i: N|_{U_i} \lm U_i \times \q, \nonumber
\end{equation}
where $\h = \p \oplus \q$.  In fact, due to the compactness of $\clm$,
for any such trivializing open cover ${\mathcal U}$, we may choose finite
refined subcovers on which we can define charts:
\begin{equation}
({\mathcal U}^j=\{U_i^j,i=1,...,s\},\sigma_i), j=1,...,J,  \ \
U^1_i \subset \overline{U_i^1}\subset
\cdot \cdot \cdot  \subset U^{J-1}_i \subset
\overline{U^{J-1}_i} \subset U^J_i \subset \overline{U^J_i}
\nonumber
\end{equation}
and $\sigma_i (U_i^j) = \b_j(0)$, the ball of radius $j$ centered about
the origin in $\p$. We simply choose for each $x \in \clm$ a
chart $(W_x,\sigma_x)$ such that $W_x \subset U_i\in {\mathcal U}$  and 
$\b_J(0) \subset \sigma_x (W_x)$ and then define $U_x^j \equiv
\sigma_x^{-1}(\b_j(0))$, and extract the finite cover. Thus, for each
$i$ and any $p \in U^j_i$, $\sigma_i(p)$ is an $m$-vector spanned
by the first $m$ elements of a fixed ordered basis of $\h$.

For the purpose of defining small local neighborhoods in the 
normal bundle, it will be convenient to consider the reduction of 
$(N,\pi|_N)$ to the Hilbert group. Recall, that $Hilb(\q)$ 
is the subgroup
of $Laut(\q)$ which preserves the inner product, so that a linear
automorphism $I$ is in $Hilb(\q)$ if and only if $I^*I=Id$ on
$\q$.  Given our Riemannian metric $g$, we may construct the Hilbert
trivializations from our vector bundle trivializations $(U_i^j,\tau
^N_i)$.  Namely, for each $x\in U_i^j$ and any $q_1$, $q_2$ in $\q$,
there exists a positive definite symmetric operator $T_{ix}$ satisfying
$g_{ix}(q_1,q_2)=\la T_{ix} q_1, q_2 \ra$.  Then, the Hilbert 
trivialization is defined for each $i$
by $(U_i^j,\tb_i)$, where $\tb_i=T_{ix}
^{\frac{1}{2}} \tau^N_{ix}$.   We have defined our Hilbert 
trivializations over the open covers $\cup_{i=1}^s U^j_i$  of $\clm$; 
thus,
any continuous map defined on some $\cup_{i=1}^s U^j_i$ can be uniformly
bounded when restricted to $\cup_{i=1}^s \overline{U^{j-1}_i}$.  
Often, we will use such compact refinements as domains of sections and 
operators.

Henceforth, $N$ will represent the reduction of the normal vector
bundle to the Hilbert group.
We define the bounded subset of the normal Hilbert bundle:
\begin{equation}
N^\eps = \{ v\in N_x: g_x(v,v) < \eps^2, \ x\in M' \}.
\nonumber
\end{equation}
$N^\eps$ defines an open neighborhood of the zero
section  of the normal Hilbert bundle;  the set of locally Lipschitz 
sections of this neighborhood $N^\eps$  will contain the 
negatively invariant manifold of the perturbed map $\sk$
upon identification of an open neighborhood of $\overline{M'}$ 
with $N^\eps$.  This identification is given by the existence of 
a tubular neighborhood diffeomorphism $h$ taking $N^\eps$, for
$\eps$ sufficiently small, onto $ V$, 
where $V$ is some $X$-open neighborhood of $\overline{M'}$.

\begin{lemma}
\label{tube}
Let $X$ be a $C^r$ manifold, $r\ge 3$, that admits a partition of
unity and let ${\mathcal M}$ be a closed submanifold.  Then  there exists
a tubular neighborhood of ${\mathcal M}$ in $X$ of class $C^{r-2}$.
\end{lemma}
This result is well known and can be found in most texts on 
infinite-dimensional manifolds (see \cite{Lang} for example).  
Of course our Hilbert manifold has partitions
of unity since the norm is differentiable; furthermore, we have a
canonical spray whose geodesic field on $T(X)$ defines the exponential
map $exp$. In other words, we require the Lagrangian vector field 
associated with the kinetic energy of the metric to exist. In general,
for infinite-dimensional Riemannian manifolds, the metric is only
weakly nondegenerate, but we will restrict attention to model
Hilbert spaces which are isomorphic to their dual space. The map $h$
is simply the restriction of the exponential to
the normal bundle, $h=exp|_N$.  
Because $h$ is at least a diffeomorphism of class $C^1$, there exists
a constant $C$ such that $\|Dh\|< C$ and $\|Dh^{-1}\| < C$ on
$\overline{V}$.  In case that $C>1$, we will use this constant in
certain global estimates, so that it may disappear  in the local
versions of these estimates.
We take our closed submanifold ${\mathcal M}$ in Lemma \ref{tube} to be 
$\overline{M'}$, and let $\sb \equiv h^{-1} \circ S \circ h$.  

We note that in many of the applications that we have in mind, the 
manifold $X$ is itself a Hilbert space, in 
which case $X$ is identified with $T(X)$ and the geodesic field
is trivial.  In that instance, the tubular neighborhood is obtained from
the smooth diffeomorphism $h:N^\eps \lm X$ defined by $(x,v) \mapsto
x+v$ and we may relax our regularity to $r \ge 2$.  In fact, if we are
content with continuous vector bundles over $\clm$, we may take 
$r\ge 1$.

\begin{defn} \label{project}
Denote by $P$ and $Q$ the following projection operators:
\begin{equation}
P:T(X)|_{M'} \lm T(M'), \text{ and } Q:T(X)|_{M'}\lm N. 
\nonumber
\end{equation}
\end{defn}

We will need the following stability condition on our negatively
invariant submanifold.
\begin{cond}[Stability]
\label{condB}
Let $B(p) \equiv Q \circ DS (S^{-1}(p))|_N:N_{S^{-1}(p)} \lm N_p$, and
let $\vartheta_1 \in (0,1)$.
Then $C^2 \|B(p)\| < \vartheta_1$ for all $p \in \cup_{i=1}^s 
\overline{U_i^{J-1}}. $
\end{cond}
This  requirement simply states that to first order, the map
$S$ decreases lengths of vectors along the normal bundle fibers
by at least a factor of $C^2\vartheta_1^{-1}$.
We show that this implies that the map $S$ sends elements of the
tubular neighborhood about $S^{-1}(M) \subset M$ ``closer'' to $M$.

\begin{lemma}
\label{lemma2}
For $0 < \eps < \eps_0$ sufficiently small, $\sb:
N^\eps|_{S^{-1}(\overline{U_i^{J-2}})} \lm N^{C^{-2}\vartheta_2 \eps}|
_{U_i^{J-1}}$ for some $\vartheta_2\in(\vartheta_1,1)$.
\end{lemma}
\paragraph{Proof.}
By definition of our tubular neighborhood, every $p\in M$ has a normal
neighborhood $R_p$ diffeomorphic to $\b_\eps(0)$ in $T_pX$. Let
$p\in \overline{U_i^{J-2}}$, and let $v\in N^\eps|_{S^{-1}(p)}$. By
continuity, we may choose $\eps$ small enough such that $h(v)\in
R_{S^{-1}(p)}$ and $S(h(v)) \in R_p$.  Since $exp_p^{-1} \circ
S \circ h$ is at least a $C^1$ map  taking $\b_\eps(0) \subset 
T_{S^{-1}(p)}X$ into $T_pX$, it differs from the linear map
$DS_{S^{-1}(p)}v$ by ${\mathcal O}(\|v\|^2)$ terms using Taylor's theorem,
and thus by ${\mathcal O}(\eps^2)$.

By projecting $exp_p^{-1} \circ S \circ h$ onto $N$ and using
Condition \ref{condB}, we may choose $\eps$ small enough so that
$C^{-2}\vartheta_1\eps+\o(\eps^2)<C^{-2}\vartheta_2\eps$, since
$C>1$.  Then, there exists some
$p' \in R_p\cap M$ such that $\|(exp|_N)^{-1}|_{p'}(S(h(v)))\|
=\|Q\  exp_p^{-1}(S(h(v))) \|$. Using the
compactness of $\overline{U_i^{J-2}}$, we may pick an $\eps_0$ small 
enough which holds for any $p\in\overline{U_i^{J-2}}$.
\quad $\blacksquare$

We will use the local Hilbert bundle trivializations that we have 
constructed to associate to each section restricted to some $U^j_i$ a
$graph$ over $\b_j(0)$ in $\p \times \q$.  This will allow us to compare
elements from distinct fibers without using parallel translation.
For each $i\in\{1,...,s\}$ and $p\in U^J_i$ and $q\in N^\eps_x$, 
define $\psi_i$ by letting
\begin{equation}
\tb_i(p,q)=(\pi|_N(p,q),\psi_i(p,q)) = (p,{\psi_i}_p(q)), 
\label{star'} 
\end{equation}
where ${\psi_i}_p(q)\in \q$.
We denote by $\s$ the sections of $N^\eps$ restricted to the cover
${\mathcal U}^{J-3}$, i.e. $\s = \N$.
Then for all $u\in \s$, we may correspond
to $u|_{U^{J-3}_i}$ the map $u_i:B_{J-3}(0)\subset \p \lm \q$ defined 
by $u_i= \psi_i \circ u \circ \sigma_i^{-1}$.  
For any two elements $v,w \in\s$, $v$ and $w$ are $\epsilon$-$C^0$ 
close if
\begin{equation}
\|v-w\|_0 :=
\max_{1\le i \le s} \sup_{x\in \b_{J-3}(0)}
\|v_i(x)-w_i(x)\| \le \epsilon.
\nonumber
\end{equation}
We define
\begin{equation}
Lip(u)=\max_{1\le i \le s} \sup_{z,z'\in \b_{J-3}(0)}
\frac{\|u_i(z)-u_i(z')\|}{\|z-z'\|}
\nonumber
\end{equation}
if it exists, and let $\sd =\{ u\in \s: Lip(u) \le \delta\}$. 
Then, for any $u\in \sd$, $h\circ Im(u)$
is a Lipschitz continuous manifold of dimension $m$.
It is important to note that in any chart overlap $U^{J-3}_i\cap U^{J-3}_k$, 
 $i,k\in \{1,...,s\}$, 
\begin{equation}
\sup_{z,z' \in \sigma_i(U^{J-3}_i)\cap\sigma_k(U^{J-3}_k)} 
\frac{\|u_i(z)-u_i(z')\|}{\|z-z'\|} < \infty
\Longrightarrow
\sup_{z,z' \in \sigma_i(U^{J-3}_i)\cap\sigma_k(U^{J-3}_k)} 
\frac{\|u_k(z)-u_k(z')\|}{\|z-z'\|} < \infty.
\nonumber
\end{equation}
Indeed, because of our reduction to the Hilbert group, there exists
a $C^{r-1}$ Hilbert automorphism $O_{ki}$ on any overlap region
$U^{J-3}_i\cap U^{J-3}_k$ so that 
\begin{equation}
\psi_k(Im(u|_{U^{J-3}_i\cap U^{J-3}_k})) = O_{ki} \ 
\psi_i(Im(u|_{U^{J-3}_i\cap U^{J-3}_k})),
\nonumber
\end{equation}
and since we were careful to shrink our patches, these Hilbert maps
are also Lipschitz in the operator norm topology.  This is a slight
generalization of the finite-dimensional setting wherein on requires
each element of an orthonormal basis to  be a $C^{r-1}$ function
over each patch, and then relies on the equivalence of norms to
conclude that the frame is $C^{r-1}$.

Finally, let $v,w \in \s$ be continuously differentiable, and to each
such section $v$ associate the $s$-vector of linear operators
$(Dv_1,...,Dv_s)$ where 
\begin{equation}
Dv_i \in C^0(\b_{J-3}(0), \l(\p,\q)) \ \ \mbox{for} \ \ i=1,...s.
\nonumber
\end{equation}
Then, $v$ and $w$ are $\epsilon$-$C^1$ close if they are 
$\epsilon$-$C^0$ close and 
\begin{equation}
\|Dv-Dw\| \equiv \max_i \sup_{x\in\b_{J-3}(0)}\|Dv_i(x)-Dw_i(x)\|
\le \epsilon. \nonumber
\end{equation}
\begin{defn}
\label{defSk}
Let $S^\kappa:X\lm X$ be $\kappa$-$C^1$ close
to $S$ on $V$, and define $\sbk \equiv h^{-1} \circ S^\kappa \circ h$.
\end{defn}

\begin{cor}
\label{cor_sk}
For $0<\eps<\eps_0$ and $\kappa > 0$ both sufficiently small,
\begin{equation}
\sbk: N^\eps|_{S^{-1}(\overline{U_i^{J-2}})} \lm N^\eps|_{U_i^{J-1}}.
\nonumber
\end{equation}
\end{cor}
\paragraph{Proof.}
This immediately follows from Definition \ref{defSk} and Lemma 
\ref{lemma2}.
\quad $\blacksquare$

Following \cite{Fen}, we will use the implicit function theorem
to show that the images of Lipschitz sections of $N^\eps$ under
$\sbk$  remain sections of the bundle.

\begin{lemma}
\label{lemma3}
For  $\delta$, $\kappa$, and $\eps$ taken sufficiently small,
$\sbk: \sd \lm \s$.
\end{lemma}
\paragraph{Proof.}
The fact that
$\si(\cup_{i=1}^sU_i^{J-1}) \subset M$ and $\si(M)\subset
\overline{M'} \subset\cup_{i=1}^s U_i^{J-3}$  together with 
Corollary \ref{cor_sk} ensures that the maps
$\chi(p;\eps, \kappa, \delta) \equiv
\pi|_N \circ \sbk \circ u \circ S^{-1}(p)$ are
well defined for all $u\in \sd$ and all $p\in \cup U^{J-2}_i$.

We must show that these maps are injective.   In any patch
$U^{J-1}_i$, we define the local representative of $\chi$ to be
$\chi_i (\eps, \kappa,\delta) = \sigma_i \circ \chi(\cdot;\eps,
\kappa,\delta)|_{U^{J-1}_i}\circ \sigma_i^{-1}$, and note that
$\chi_i(0, 0,\delta)$ is the inclusion map of $\b_{J-2}(0)$ into $\p$.
By the implicit function theorem for Lipschitz continuous maps (see
\cite{Die}), there exists a 
neighborhood ${\mathcal X}$ of $\chi_i(0,0,\delta)$ in the 
Lipschitz topology  such that each $f \in {\mathcal X}$ maps 
$\b_{J-2}(0)$ into $\p$  injectively and satisfies
\begin{equation}
\overline{\b_{J-3}(0)}\subset f(\b_{J-2}(0)) \subset 
\overline{f(\b_{J-2}(0))} \subset \b_{J-1}(0).
\nonumber
\end{equation}
To each $f \in {\mathcal X}$ we identify the injective map $\bar{f}$ on
$U_i^{J-1}$ and have 
\begin{equation}
\overline{U_i^{J-3}} \subset \bar{f}(U_i^{J-2}) \subset 
\overline{\bar{f}(U_i^{J-2})} \subset U_i^{J-1}.
\label{star}
\end{equation}
By taking $\eps$, $\delta$, and $\kappa$ small enough, 
$\chi_i(\eps,\kappa,
\delta)$ can be made arbitrarily close to the inclusion map
and thus is an element of ${\mathcal X}$.  Since this is true for each
$i$, we see that
each $\chi(\cdot;\eps,\kappa,\delta)$ is injective and that
$\overline{\cup U_i^{J-3}} \subset \chi (\cup U^{J-2}_i)$.  It is for
this reason that we have defined our sections on $\cup U^{J-3}_i$.
\quad $\blacksquare$

In order to show that $\sk$ has an overflowing invariant manifold,
we will prove that the map $\sbk$ takes $\sd$ into itself  (more
precisely, $\sbk$ maps the images of elements of $\sd$)
and is in fact a contraction.  To do so, we will need to consider the
partial derivatives of its local representation.  For each point
$(p,v) \in N^\eps|_{U_i^J}$, let $L_i:N^\eps|_{U_i^J} \lm \p \times
\q$ be defined by 
\begin{equation}
L_i(p,v)=(\sigma_i(p),\psi_i(p,v)) \equiv (x,y). \nonumber
\end{equation}
Then $L_i^{-1}((x,0))$ is an element of the zero section of $N^\eps$.
By compactness,  for any $i\in \{ 1,...,s\}$, we may choose the
constant $C$ in Condition \ref{condB} large enough so that we 
have the uniform bounds
$\|D\sigma_i\| < C$ and $\|D\sigma_i^{-1}\| < C$
on $\cup \overline{U_i^{J-1}}$, and by definition of the Hilbert
group, $\|D\psi_i (p,v)\| = \|D\psi^{-1}_i(\psi_i(p,v))\| = 1$ for all
$(p,v) \in N^\eps|_{\cup U_i^J}$. 

\begin{cond}[Smoothing]
\label{condA}
Let $A(p) = DS^{-1}|_{M'}(p):T_pM' \lm T_{S^{-1}(p)}M'$, and let
$\vartheta_1\in (0,1)$.  Then
$C^{4} \|A(p)\| \|B(p)\| < \vartheta_1$ for all $p \in \cup
\overline{U_i^{J-1}}$.
\end{cond}

In the case that $X$ is locally compact, $DS$ is uniformly bounded
and $D(S|_M')$ has uniformly bounded inverse on bounded neighborhoods
of $M'$.  In the more general case, we may choose a sufficiently
small neighborhood of $M'$ contained in the tubular neighborhood on
which the same conclusion holds.

\begin{lemma}[Invertability]
\label{invert}
There exists a sufficiently small $X$-open neighborhood $V$ of 
$\overline{M'}$ on which $DS$ is uniformly bounded and
$D(S|_M)$ has uniformly bounded inverse.
\end{lemma}
\paragraph{Proof.}
Since $S$ is $C^1$ and $\overline{M'}$ is compact, there exists
$C>0$ such that $\|DS\| < C$ on $\overline{M'}$.  For any $\delta >0$
and each $p\in \overline{M'}$, we may choose $\eps(p)$ such
that $\|DS \| < C + \delta$ on ${\mathcal B}_{\eps(p)}(p)$.  Set 
$V= \cup_{p\in {\mathcal P}} {\mathcal B}_\eps (p)$, where ${\mathcal P}$
is a finite subset of points on $\overline{M'}$ provided by compactness,
and $\eps = \min_{p \in {\cap P}} \eps(p)$.

Our conditions ensure that $S|_{\overline{M'}}$ is a diffeomorphism,
and since
$D(S|_{\overline{M'}}^{-1})$ is uniformly bounded on $\overline{M'}$,
$[D(S|_{\overline{M'}})]^{-1}$ exists and is uniformly bounded by
the inverse function theorem.  Thus, by the uniform continuity
of the spectrum of $DS$ and the compactness of $\overline{M'}$,
we may shrink $V$ if necessary so that 
$[D(S|_{\overline{M'}})]^{-1}$ is uniformly bounded on $V$.
\quad $\blacksquare$

For each $(x,y) \in \p \times \q$ satisfying $L_i^{-1}(x,y) \in
N^\eps|_{\overline{U_i^{J-2}}\cap S^{-1}(\overline{U_k^{J-2}}) }$,
we can define the local {\it graph transform} of a section of $N^\eps$,
represented locally as a graph over $\overline{U_i^{J-2}}$, into a
graph over $\overline{U_k^{J-2}}$ by $(x,y) \mapsto \G (x,y)$ where
$\G : \p\times \q \lm \p\times\q$ and is defined by
\begin{equation}
\G = (\Gx , \Gy )\equiv (\sigma_k \circ \pi|_N \circ \sbk\circ L_i^{-1},
\psi_k \circ \sbk \circ L_i^{-1}).
\label{local}
\end{equation}
Using this notation, we may translate the result of Corollary
\ref{cor_sk} into the local form
\begin{equation}
\| \Gy (x,y) \| < \eps , \label{1}
\end{equation}
and we may use Definition \ref{defSk} together with Conditions 
\ref{condB} and
\ref{condA} to get for $\kappa$ sufficiently small
\begin{equation}
\| D_2 \Gy (x,y)\| < \vartheta_2, \ \ \ \ \
\| (D_1 \Gx (x,y))^{-1} \| \|D_2 \Gy (x,y)\| < \vartheta_2, \ \ \
\vartheta_2 \in (\vartheta_1,1).
\label{2}
\end{equation}
Since  $M'$ is negatively invariant, 
${D_1 G_{ik}^0}^2(x,0)= 0$, so for $\eta>0$, we may choose
$\eps$ and $\kappa$ small enough such that
\begin{equation}
\|D_1 \Gy(x,y)\| < \eta . \label{3}
\end{equation}
Finally, Lemma \ref{invert} gives us
\begin{equation}
\label{4}
\| D \G (x,y) \| < \tilde{C}, \ \mbox{and} \ 
\|[D_1\Gx (x,y)]^{-1}\| < \tilde{C} \ \ \forall (x,y) \in
L_i(N^\eps|_{\overline{U_i^{J-2}}\cap S^{-1}(\overline{U_k^{J-2}}) }).
\end{equation}
for some bounded constant $\tilde{C}$.
The first inequality is valid for a sufficiently small $\kappa$ because
of the $C^1$ closeness of $\sk$ with $S$ assumed in Definition 
\ref{defSk}.  For the second inequality, we use the uniform continuity 
of the
spectrum; namely, if the spectrum of $D_1 {G_{ik}^0}^1(x,y)$ does not
intersect a neighborhood of zero as given by Lemma \ref{invert}, 
then we may choose $\kappa$ small enough so that the same is true 
for $D_1 \Gx (x,y)$.  

In the case that a section
$u\in\sd$  is continuously differentiable, we have that
\begin{equation}
\label{4a}
\|[D\Gx (x,u_i(x))]^{-1}\| < \tilde{C} .
\end{equation}

We make a remark on the constant $C$.
Our canonical spray defines a unique bilinear form which in
turn gives rise to a unique Riemannian connection.  For simplicity of
analysis, we have pulled back this connection in each local 
trivialization by isomorphisms which we have bounded by $C$.

\subsection{Lipschitz negatively invariant manifolds for the perturbed
mapping}

The proof of the following lemma relies heavily on the smoothing
condition \ref{condA} in its local form (\ref{2}) which is essential for
overcoming the possibly large bound $\tilde{C}$ on the norm of the 
local derivatives of $\G$ in (\ref{4}).
\begin{lemma}
\label{Lip}
For $\eps$, $\delta$, and $\kappa$ sufficiently small, $\sbk :\sd \lm
\sd$.
\end{lemma}
\paragraph{Proof.}
To prove that $\sbk (\sd ) \subset \sd$, it is equivalent to use the
local representation (\ref{local}), and show that for all $u \in \sd$
\begin{equation}
\| \Gy (x, u_i(x)) - \Gy (x',u_i(x')) \| \le \delta
\| \Gx (x, u_i(x)) - \Gx (x',u_i(x')) \|
\label{6}
\end{equation}
for $x, x' \in L(N^\eps_{U_i^{J-3} \cap S^{-1}(U_k^{J-3})})$.
As usual, we will get a lower bound on the right hand side of 
equation (\ref{6}) in terms of the difference of $x$ and $x'$.
We have that
\begin{eqnarray*}
\| \Gx (x,u_i(x)) - \Gx (x',u_i(x'))\|  &  \ge  &
\| \Gx (x,u_i(x)) - \Gx (x',u_i(x))\| \\ && - \, \, 
\| \Gx (x',u_i(x)) - \Gx (x',u_i(x'))\|  \\ & \equiv &
\| \t_1 \| - \|\t_2\|.
\end{eqnarray*}
With Lemma \ref{invert} in its local form (\ref{4}),
$D_1\Gx$ is invertible, and so we may 
use the implicit function theorem together with a Taylor expansion to 
get a lower bound on $\t_1$; namely, since
\begin{equation}
x-x' = [D_1\Gx]^{-1} (\Gx (x,u_i(x)) - \Gx (x',u_i(x)) + \o(\|x-x'\|^2),
\nonumber
\end{equation}
we have that
\begin{equation}
\|\t_1\| \ge \| [D_1\Gx]^{-1}\|^{-1} \|x -x'\| + \o(\|x-x'\|^2).
\label{7}
\end{equation}
The upper bound on $\|\t_2\|$ is simply
\begin{equation}
\| \t_2\| \le \delta \|D_2 \Gx \| \| x-x'\| + \o(\|x-x'\|^2).
\nonumber
\end{equation}
We subtract $\| \t_2\|$ from $\|\t_1\|$, factor the term 
$\| [D_1\Gx]^{-1}\|^{-1}$, use the bounds in (\ref{4}), and take
$x'$ close enough to $x$ so that $\o(\|x - x'\|) < \delta \tilde{C}
\|x - x'\|$. Then 
\begin{equation}
\|\t_1\| - \|\t_2\| \ge
\| [D_1 \Gx]^{-1} \|^{-1} (1 - 2 \delta \tilde{C}^2) \|x - x'\|.
\label{star''}
\end{equation}
Let $\vartheta_2 < \varrho < 1$, and choose $\delta < \frac{1-\varrho}
{2\tilde{C}^2}$.  This ensures that the constant multiplying 
$\| x - x' \|$ is strictly positive.

Next, we estimate
\begin{eqnarray*}
\| \Gy (x,u_i(x)) - \Gy (x', u_i(x')) \| & \le &
(\|  D_1 \Gy \| + \delta \|D_2 \Gy \| ) \|x -x'\| + \o (\|x-x'\|^2).
\end{eqnarray*}
Then, we take $\kappa$ small enough so that (\ref{3}) holds 
with $\eta < \frac{\varrho-\vartheta_2}{2\tilde{C}}$ and
choose $\| x -x'\|$ even smaller if necessary so that 
$\o (\| x - x' \|^2) \le  \eta \|x -x'\|$. We get
\begin{eqnarray*}
\| \Gy (x,u_i(x)) - \Gy (x', u_i(x')) \| & \le &
\frac{1}{\varrho}(2 \eta \|[D_1\Gx]^{-1} \|^{-1}
+ \delta\|D_2\Gy\| \  \|[D_1\Gx]^{-1} \|^{-1}) 
\\ && \cdot  \| \Gx (x,u_i(x))
- \Gx (x',u_i(x')) \|,
\end{eqnarray*}
and using (\ref{2}) and (\ref{4}) completes the estimate.  Finally,
we note, as in \cite{Fen}, that since $\overline{\b_{J-3}(0)}$ is 
compact and convex, the estimate we have just derived holds in the
entire ball.  Thus, from our definition of $\sd$, we have shown that
$\sbk (u)$ is Lipschitz.
\quad $\blacksquare$
\begin{thm}
\label{thm1}
Assume Conditions \ref{condB} and \ref{condA} hold,
and let $\eps$, $\delta$, and
$\kappa$ be sufficiently small. Then there exist Lipschitz submanifolds 
$\overline{M_\kappa}$  of $X$ which are negatively invariant with
respect to the perturbed mappings $S^\kappa$. Furthermore, 
$\overline{M_\kappa} \lm \clm$ as $\kappa \lm 0$ in $C^0$.
\end{thm}
\paragraph{Proof.}
The proof proceeds in the usual way.  We first show that  $\sbk$ is
a contraction  on $\sd$ in the $C^0$ topology.  Then, since $\sd$
is a metric space, closed under $C^0$ convergence, we appeal to the
contraction mapping  theorem to show that $\sbk$ has a fixed
point $u^\kappa\in\sd$,  and define $\overline{M_\kappa} =
h(Im(u^\kappa|_{\clm}))$.

Thus, we will obtain a uniform bound on the distance between the image
of two sections along any fiber over $\cup U_i^{J-3}$.  Let $u, u'
\in \sd$;  then, for any $\bar{p}$ in some $U_k^{J-3}$, the proof of
Lemma \ref{lemma3} allows us to choose some $i$ and $x,x'$  in 
$\sigma_i(U_i^{J-3})$ such that 
\begin{equation}
\label{8}
\sigma_i(\bar{p}) =\Gx (x,u_i(x))= \Gx(x',u'_i(x')).
\end{equation}
A simple estimate shows that
\begin{eqnarray}
\| \Gy (x,u_i(x)) -\Gy (x',u'_i(x'))\| & \le & 
 \|D_1\Gy\|\cdot\|x-x'\| + \o(\|x-x'\|^2) \nonumber \\
&& \, \, + \|D_2\Gy\|\cdot\|u_i(x)-u'_i(x)\|+\o(\|u_i(x)-u'_i(x)\|^2)
\nonumber \\
&& \, \, + \|D_2\Gy\|\cdot\|u'_i(x)-u'_i(x')\|+
           \o(\|u'_i(x)-u'_i(x')\|^2) \nonumber\\
&\le& (2\eta + \vartheta_2 \delta)\|x - x'\| + \vartheta_2
\|u_i - u'_i\|_0 + \o(\|u_i-u'_i\|_0^2). \label{8a}
\end{eqnarray}
Taking $x'$ close enough to $x$ so that $\o(\|x-x'\|^2) < 
\frac{1}{2} \| [D_1\Gx ]^{-1}\|^{-1}\|x-x'\|$ and using (\ref{7}), we
have
\begin{equation}
\| \t_1\| \ge \frac{1}{2} \| [D_1\Gx ]^{-1}\|^{-1}\|x-x'\|.
\label{9}
\end{equation}
By (\ref{8}), $\| \t_1\|=\|\Gx (x',u'_i(x')) - \Gx (x',u_i(x)) \|$ so
we get
\begin{eqnarray*}
\frac{1}{2} \| [D_1\Gx ]^{-1}\|^{-1}\|x-x'\| &\le & \|D_2 \Gx\|
\| u'_i(x')-u'_i(x) \| + \o (\|u'_i(x')-u'_i(x)\|^2) \\
&& \, \, +\|D_2 \Gx\| \|u'_i(x) -u_i(x)\| + \o(\|u'_i(x) -u_i(x)\|^2) \\
&\le& \tilde{C} (\delta \|x -x'\| + \|u'-u\|_0) + \o(\|x-x'\|^2)
+\o(\|u'_i - u_i\|_0^2).
\end{eqnarray*}
Then, since $u,u'\in\sd$, (\ref{8a}) permits us to choose $\delta$ and
$\kappa$ small enough and $\o(\|u - u'\|_0^2) < (1-\vartheta_2)
\|u_i-u'_i\|_0$ to get the desired contraction.

Finally, since locally, $u^0_i$ and $u^\kappa_i$ are the fixed points
of $G^0_{ik}$ and $G^\kappa_{ik}$, respectively, we may use the contraction
property of these maps (which we have just proven) to show that
$u^\kappa \lm u^0$ in $C^0$ and hence that $\overline{M_\kappa}\lm \clm$ in
$C^0$.
\quad $\blacksquare$

\subsection{$C^1$ negatively invariant manifold for  the perturbed 
mapping}
If $u \in \sd$ is continuously differentiable and $p \in 
\cup_{i=1}^sU_i^{J-3}$, then 
$Du(p)\in \l (T_pM', T_{u(p)}N)$ can locally be represented 
in each $U_i^{J-3}$ by
\begin{equation}
Du_i \in C^0(\b_{J-3}(0), \l (\p,\q)) \ \ \ \mbox{for} \ \ \ i=1,...s.
\nonumber
\end{equation}
Hence, we associate to $Du$ the $s$-vector $(Du_1,...,Du_s) \in
[ C^0(\b_{J-3}(0), \l (\p,\q))]^s $  which is a complete metric space
when normed by 
\begin{equation}
\|Du\| \equiv \max_i \sup_{x\in \b_{J-3}(0)} \|Du_i(x)\|,
\label{norm}
\end{equation}
having the topology of uniform convergence.  We then define the
subset 
\begin{equation}
\gd = \{T \in [ C^0(\b_{J-3}(0), \l (\p,\q))]^s:\| T\| \le \delta\}.
\nonumber
\end{equation}

We let $u^\kappa$ be the Lipschitz section in $\sd$ which defines
the negatively invariant manifold of the perturbed mapping $\sk$.
The invariance can locally be represented as
\begin{equation}
u^\kappa_k(\Gx(x,u^\kappa_i(x)))= \Gy(x,u^\kappa_i(x)). \label{10}
\end{equation}
In order to show that $\uk$ is actually of class $C^1$, we will 
construct  a Cauchy sequence in $\sd$ which converges to a fixed
point $\uk \in \sd$, and then prove that the sequence is in fact
Cauchy in $C^1$. In the case that $\uk$ is Fr\'{e}chet differentiable,
we may differentiate (\ref{10}), and find that its $s$ components 
satisfy
\begin{eqnarray}
D\uk_k (\Gx(x,\uk_i(x)))& =& \left[ D_1\Gy(x,\uk_i)) + D_2 \Gy (x,\uk_i(x))
D\uk_i(x) \right] \nonumber \\
&& \, \, \cdot \ \left[ D_1\Gx(x,\uk_i)) + D_2 \Gx (x,\uk_i(x))
D\uk_i(x) \right]^{-1} \label{Du} \\
&\equiv& \left[ \lki^{\uk} (D\uk_i)\right] \left\{ (\Gx(x,\uk_i(x))
\right\}, \nonumber
\end{eqnarray}
where the superscript $\uk$ on $\lki$ refers to the section on which
the partial derivatives of $\Gx$ and $\Gy$ are evaluated.  We note
that (\ref{Du}) is well defined because of Lemma \ref{invert} in
its local form (\ref{4a}).  Thus, in the case that $\uk$ is $C^1$, we
have computed the derivative of the local graph transform, which maps
elements $(x,y) \in L_i(N^\eps|_{\overline{U_i^{J-2}}\cap 
S^{-1}(\overline{U_k^{J-2}}) })$  into $\b_{J-1}(0)$. In order to obtain
a global representation for the derivative of $\uk$ we must ``piece''
together our local representations using a partition of unity argument.
This, in effect, allows us to use our local vector space structure to
define a global derivative, without having to analyze the spray 
associated to the true covariant derivative.  We will need the 
following:
\begin{defn}
\label{unity}
A partition of unity of class $C^r$ on a manifold $\m$ consists of an
open covering $\{U_i\}$ of $\m$ and a family of functions $\phi_i\in
C^r( \m,\r )$ satisfying the following conditions:
\begin{itemize}
\item[1.] for all $m\in \m$, $\phi_i(m)\ge 0$,
\item[2.] the support of $\phi_i$ is contained in $U_i$,
\item[3.] the covering is locally finite, and
\item[4.] for each point $m\in\m$, $\sum \phi_i(m)=1$.
\end{itemize}
\end{defn}

This will be used to unite the $s$ local representations of the
derivative into a single operator equation.
By Urysohn's lemma, for each $p$ in a $C^r$ manifold $\m$ and
$W \subset \m$ a neighborhood of $p$ contained in a coordinate
neighborhood  of $p$ which is diffeomorphic to an open ball, there 
exists a neighborhood $U$ of $p$ with $\overline{U}\subset W$ and
a $C^r$ function $\phi:\m \lm \r$ with $0\le \phi(m) \le 1$ if $m\in\m$,
$\phi(m) = 1$ if $m\in \overline{U}$, and $\phi(m)=0$ if $m\in W^c$.
Thus, choose $\gamma_i=1$ on $U_i^{J-5}$ with $supp(\gamma_i)\subset
U_i^{J-2}$ and since the covering $\{U_i^{J-4}\}$ is locally finite,
so is $\{ \gamma_i \}$; then,
simply let $\phi_i(p) = \gamma_i(p)/ \sum_{k=1}
^s \gamma_k(p)$ to satisfy Definition \ref{unity} with $\m = \cup U_i^
{J-3}$.  From (\ref{Du}) it is clear that 
\begin{equation}
D\uk_k (x) = \sum_{i=1}^s \phi_i(p) \lki^{\uk}(D\uk_i)(x), \ \ k=1,...s,
\label{derivative}
\end{equation}
where $p\in M$ satisfies $\sbk(\uk(p)) = \uk(\sigma_k^{-1}(x))$.
The patched-together derivative displayed in (\ref{derivative}) serves 
as motivation for our next two lemmas.

To simplify our notation for a fixed section $u$, we shall define
$d_rg^s(\xi) = D_r{\G}^s(\xi, u_i(\xi))$ for $r,s=1,2$.
\begin{lemma}
\label{lemma5}
Let $T^1=(T^1_1,...,T^1_s) \in \gd$, choose $u \in \sd$,  and let 
$T^2$ be defined by
\begin{equation}
T^2_k (x) = \sum_{i=1}^s \phi_i(p) \lki^u (T^1_i)(x),
\label{pound}
\end{equation}
where for each $p \in U^{J-3}_i$, $x=\Gx(\sigma_i(p),u(\sigma_i(p)))$.
Then for all $\eps$, $\delta$, and
$\kappa$ sufficiently small, $T^2 \in \gd$.
\end{lemma}
(We note that since $\| T^1 \| < \delta$, (\ref{pound}) is well defined 
because of Lemma \ref{invert} in its local form (\ref{4a}).) 
\paragraph{Proof.} 
We first show that $\|T^2\|\le \delta$, and to do so, we use  
the estimates (\ref{2})-(\ref{4}).
Since the norm on $\gd$ in (\ref{norm}) is defined by computing
the maximum over $i\in\{1,...,s\}$, 
it is sufficient to show 
that for all $i,k \in \{ 1,...,s\}$, $\| \lki^u(T^1_i)(x)\| \le
\delta$ for all $x\in U_i^{J-3}$.
We have
\begin{eqnarray*}
T^2_k & = & \lki^u (T^1_i) \\
&=& \left[ d_1g^2 + d_2g^2 T^1_i\right] \left[d_1g^1 + d_2g^1 T^1_i
\right]^{-1} \\
&=& \left[ d_1g^2 + d_2g^2 T^1_i\right] \left[Id + d_2g^1 T_i^1 
(d_1g^1)^{-1} \right]^{-1} (d_1g^1)^{-1},
\end{eqnarray*}
so $\| T^2_k\| \le \|d_1g^2 + d_2g^2 T^1_i \| \cdot
\|(d_1g^1)^{-1}\| \cdot \| id + d_2g^1 T^1_i (d_1g^1)^{-1}\|$.
Also, 
\begin{eqnarray*}
\| d_1g^2 + d_2 g^2 T^1_i \|  \cdot \|(d_1g^1)^{-1}\|
&\le& \left\{ \|d_1g^2\| + \|d_2g^2\| \cdot \| T_i^1\| \right\}
\cdot \|(d_1g^1)\|^{-1}\\
&\le& \tilde{C} \eta + \vartheta_2 \delta,
\end{eqnarray*}
and
\begin{equation}
\| Id + d_2g^1 \ T^1_i (d_1g^1)^{-1} \| < 1 + \tilde{C}^2 \delta
+\o(\delta^2)
\nonumber
\end{equation}
since $\| d_2g^1 \cdot T^1_i \cdot (d_1g^1)^{-1}\| < 1$ for 
$\delta$ small enough. Therefore, 
\begin{equation}
\|T^2_k\| \le \eta(\tilde{C}+\tilde{C}^3 \delta)
+\vartheta_2 \delta + \o(\delta^2) + \o(\eta \delta) \le \delta
\nonumber
\end{equation}
for $\kappa$ and hence $\eta$ sufficiently small.

Next, we must show that $T^2_k$ is continuous.  The domain of $\Gx$ is
$N^\eps|_{\overline{U}_i^{J-2}\cap S^{-1}(\overline{U}_k^{J-2})}$.
Choose $k\in\{1,...,s\}$ such that for $\xi\in\b_{J-3}(0)$ and
$\sigma_k^{-1}(\xi)\in U_k^{J-3}$, $\xi = \Gx(x,u_i(x))$ for some
$i\in\{1,...,s\}$, and $x\in \b_{J-3}(0)$ and $\sigma_k^{-1}(\xi)\in
\overline{U}_k^{J-1}$.  We may choose such a $k$, because of the 
nesting property of the $U^j_k$'s and (\ref{star}).  Choose the largest
$c>0$ so that $\b_c(\xi) \subset \b_{J-3}(0)$.  By compactness, we may
choose $c$ independent of $i,k$ or $\xi$.  

Then for any $d>0$, we must show that for 
\begin{equation}
\| \xi - \xi'\| < c \ \ \mbox{taken sufficiently small},
\label{21}
\end{equation}
$\|T^2_k(\xi)-T^2_k(\xi')\| <d$, where $\xi'=\Gx(x',u_i(x'))$.

This result relies on the continuity of $T^1_i$.  Let
\begin{equation}
A^s(x) = d_1g^s (x) + d_sg^s(x) T^1_i(x), \ \ \ s=1,2.
\nonumber
\end{equation}
Then
\begin{eqnarray*}
\|T^2_k(\xi)-T^2_k(\xi')\| &\le& \|A^2(x)\|\cdot \|A^1(x)^{-1}-A^1(x')
^{-1}\|  + \|A^2(x)-A^2(x')\| \cdot \|A^1(x')^{-1}\|\\
&\le& \|A^2(x)\|\cdot\|A^1(x')^{-1}\|\cdot \|A^1(x)-A^1(x')\|\cdot
\|A^1(x')^{-1}\| \\
&&\, \, + \|A^1(x')^{-1}\|\cdot \|A^2(x)-A^2(x')\|,
\end{eqnarray*}
where ${A^1}^{-1}$ exists due to Lemma \ref{invert}.

Using (\ref{star''}), we see that (\ref{21}) implies that $\|x-x'\|$
can be made arbitrarily small.  Since $\|A^2\|$ and $\|{A^1}^{-1}\|$
are bounded due to (\ref{4}), the continuity of $d_rg^s$ for $r,s=1,2$, 
and $T^1_i$ give us the result.
\quad $\blacksquare$
\begin{lemma}
\label{lemma6}
(a) The mapping $\sum \phi_i \lki^u$ is a contraction on $\gd$ for each
$u\in\sd$ with a unique fixed point $T^u$, and 
(b) the map $u\mapsto T^u$ is continuous.
\end{lemma}
\paragraph{Proof.}
(a) Fix $u\in \sd$ and let $T^u$ and $S^u$ be in $\gd$.  Then
\begin{eqnarray*}
\lefteqn{\left\|\lki^u (T^u_k) - \lki^u(S^u_k)\right\|}\\
&& = \left\| [d_1g^2 + d_2g^2 T_i^u][d_1g^1+d_2g^1 T^u_i]^{-1} -
[d_1g^2 + d_2g^2 S_i^u][d_1g^1+d_2g^1 S^u_i]^{-1} \right\| \\
&& =\left\| [d_1g^2 + d_2g^2T^u_i][d_1g^1+d_2g^2 T_i^u]^{-1}
\left\{[d_1g^1 + d_2g^1S^u_i] - [d_1g^1 + d_2g^1T_i^u]\right\}
[d_1g^1 + d_2g^1 S^u_i]^{-1} \right.\\
&&\,\, +\left. \left\{ [d_1g^2 T_i^u] - [d_1g^2 + d_2g^2 S_i^u] \right\}
[d_1g^1 + d_2g^1 S^u_i]^{-1}\right\|.
\end{eqnarray*}
Using the same type of estimate as in Lemma \ref{lemma5}, we take
$\kappa$ small enough for a sufficiently small $\eta$  to get
\begin{eqnarray*}
\left\| \lki^u (T^u_k) - \lki^u(S^u_k)\right\| &\le&
\left( \| d_2g^2 \| \ \|(d_1g^1)^{-1} \| + \delta \| d_2g^1 \| 
\ \|(d_1g^1)^{-1}\| \right)  \{ 1+ \o(\delta) \} \|T_i^u-S_i^u\| \\
&\le& (\vartheta_2 + \o(\delta)) \| T^u_i - S^u_i\| \le \mu 
\|T^u_i - S_i^u\|, \ \ \ \mu < 1,
\end{eqnarray*}
when we take $\delta$ small enough for $\o(\delta) \le 1 - \vartheta_2$,
and $\mu$ is the maximum of all the contraction coefficients obtained
in each chart.

(b) For $r=1,2$, let $T^{u^r}$ be the fixed point of $\sum \phi_i
\lki^{u^r}$.  Choose $u^1$ and $u^2$ in $\sd$ and for all $n\in\n$, let
the sequences $\{ \tbn\}$ and $\{\tbbn\}$ be defined by
\begin{equation}
\overline{T}^{n+1}_k = \sum \phi_i \lki^{u^1} (\tbni), \ \ \
\overline{\overline{T}}^{n+1}_k = \sum \phi_i \lki^{u^2} (\tbbni),
\nonumber
\end{equation}
where $\overline{T}_k^0$ and $\overline{\overline{T}}^0_k$ are chosen
such that $\| \overline{T}^0 - \overline{\overline{T}}^0\| =
\o (\|u^1 - u^2\|_0)$.   Using part (a), we see that the sequences are
Cauchy and must have limits $T^{u^1}$ and $T^{u^2}$, respectively.
Thus, in order to prove that $\| T^{u^1}- T^{u^2}\|=\o(\|u^1-u^2\|_o)$,
we may use induction, and show that if $\| \overline{T}_k^{n-1} -
\overline{\overline{T}}_k^{n-1}\| = \o(\|u^1-u^2\|_o)$, then
$\| \tbnk - \tbbnk\|$ remains $\o(\|u^1-u^2\|_o)$.  The proof of this
estimate is similar to those already shown above.
\quad $\blacksquare$

\begin{thm}
\label{thm2}
Assume Conditions \ref{condB} and \ref{condA} hold,
and let $\eps$, $\delta$, and
$\kappa$ be sufficiently small. Then the Lipschitz negatively invariant
manifolds $\overline{M_\kappa}$ of the perturbed mappings $\sk$ are
$C^1$, and $\overline{M_\kappa} \lm \clm$ as $\kappa \lm 0$ in $C^1$.
\end{thm}
\paragraph{Proof.}
Let $\uk$ be the Lipschitz section supplied by Theorem \ref{thm1}
such that $h \circ \uk |_{\clm} = \overline{M_\kappa}$, and let the
components of $T^{\uk}$ satisfy the functional equation $T^{\uk}_k
= \sum \phi_i \lki^{\uk} (T^{\uk})$.  We will show that $\uk$ is 
continuously differentiable and that $T^{\uk} = D\uk$.

We define the sequence $\{u^n\} \subset \sd$ iteratively by
$u^{n+1} = \sbk (u^n)$ so that locally
\begin{equation}
u^{n+1}_k (\Gx (x,u^n_i(x)) = \Gy(x,u^n_i(x)). \label{22}
\end{equation}
If $u^n\in \sd$ and is of class $C^1$, then $u^{n+1} \in \sd$ is also
$C^1$ using (\ref{4a}) together with an implicit function theorem
argument.  Thus, if we choose $u^0 \in \sd$ such that $Du^0 \in \gd$,
($u^0=0$, for example), then for all $n\in\n$,
$Du^n \in
[C^0(\b_{J-3}(0),\l(\p,\q))]^s$.  By computing the derivative
of (\ref{22}), we see that
\begin{equation}
Du_k^{n+1}= \sum_{i=1}^s \phi_i \lki^{u^n} (Du^n_i). \nonumber
\end{equation} From the proof of Lemma \ref{lemma5}, it is clear 
that $Du^n \in\gd$
for all $n\in\n$.  Theorem \ref{thm1} shows that $\sbk$ is a contraction
on $\sd$ so $\{u^n\}$ is Cauchy in $\sd$ and has as its unique limit
$\uk$; therefore, to prove that $\{u^n\}$ is Cauchy in $C^1$, it 
suffices to show that $\{Du^n\}$ is Cauchy in $\gd$.

We argue as in [\cite{Chow}, Lemma 4.1].  From part (a) of Lemma
\ref{lemma6}, we may choose, for each $n$, $T^{u^n} \in \gd$ satisfying
\begin{equation}
T^{u^n}_k = \sum_{i=1}^s \phi_i \lki^{u^n} (T_i^{u^n}), \nonumber
\end{equation}
and define $E_N = sup_{n,m>N}\|T^{u^n}-T^{u^m}\|$.
Since $\sup_{n,m > N} \|u^n - u^m\|_0 \lm 0$ as $N\lm\infty$, Lemma
\ref{lemma6}(b) states that $E_N \lm 0$ as $N\lm \infty$.
Now, for each $k\in \{1,...,s\}$,
\begin{equation}
\| Du^n_k - Du^m_k\| \le \|Du^n_k - T^{u^n}_k\| +
\| T^{u^n}_k - T^{u^m}_k\| +\| T^{u^m}_k - Du^m_k\|,  \label{triangle}
\end{equation}
and
\begin{eqnarray*}
\|Du^{n+1}_k - T^{u^{n+1}}_k\| &\le& \|Du^{n+1}_k - T^{u^n}_k\| + E_N\\
&\le& \| \sum \phi_i \lki^{u^n}(Du^n_i) - 
\sum \phi_i\lki^{u^n}(T^{u^n})\| + E_N\\
&\le& \mu \|Du^n_i-T^{u^n}_i\| + E_N,
\end{eqnarray*}
where $\mu$ is the contraction coefficient supplied by Lemma
\ref{lemma6}(a).

Proceeding by induction, we find that for $n\ge m \ge N$,
\begin{equation}
\| Du^n - T^{u^n}\| \le \mu^{m-N} \|Du^N -T^{u^N}\| + \frac{E_N}{1-\mu},
\label{EN}
\end{equation}
so using (\ref{triangle}), for $n\ge m \ge N$,
\begin{equation}
\|Du^n - Du^m \| \le 2\mu^{m-N} \|Du^N - T^{u^N}\| +
2\frac{E_N}{1-\mu} + E_N. \label{EN1}
\end{equation}
Then, for any $\eps >0$, we may choose $N$ large enough so that 
each of the
last two term on the right hand side of (\ref{EN1}) are less than
$\eps/3$.  We then choose $m>N$ large enough so that the first term
on the right hand side of (\ref{EN1}) is less than $\eps/3$.
This shows that the sequence $\{u^n\}$ is $C^1$-Cauchy.  We 
conclude, using the completeness of $C^1$,
that $u^n \lm \uk $ in $\sd$ and $Du^n \lm D\uk$ in $\gd$ as
$n\lm \infty$, and by passing to the limit in (\ref{EN}) we obtain 
$D\uk = T^{\uk}$.

Finally, the contraction property of part (a) of Lemma 
\ref{lemma6} assures us that $\|D\uk -Du^0\|\lm 0$ as $\kappa\lm 0$,
and together with Theorem \ref{thm1}, we obtain $\overline{M_\kappa} 
\lm \clm$ in $C^1$.
\quad $\blacksquare$

Next, assume our mapping $S$ arises as the time-$T$ map of a nonlinear
$C^1$ semi-group $S_t$.  
As in \cite{Fen}, we define the following generalized
Lyapunov-type numbers
\begin{equation}
\begin{array}{c}
\kk(x) = \limsup_{t\lm\infty} \|B_t(x)\|^{1/t} \\
\theta(x) = \limsup_{t\lm\infty} \frac{\log \|A_t(x)\| }
                                      {-\log \|B_t(x)\|},
\ \ if \ \ \kk(x) < 1,
\end{array}
\nonumber
\end{equation}
where $A_t(x) \equiv DS_{-t}|_{M'}(x):T_xM' \lm T_{S_{-t}(x)}
M'$ and $B_t(x)\equiv QDS_t(S_{-t}(x))|_N: N_{S_{-t}}(x) \lm
N_x$.  We have that $\kk$ and $\theta$ are bounded on $\clm$.  
Furthermore, it is proven in
\cite{Fen} that the generalized Lyapunov-type numbers are constant on
orbits and hence depend only on the backward limit sets on $M$.  In
addition, they give us uniform estimates on the norms of $A_t(\cdot)$
and $B_t(\cdot)$ as was shown by Sacker and later Fenichel (see
\cite{Sac}, \cite{Fen}) in the following uniformity lemma.

\begin{lemma}
\label{uniformity}
a) If $\lim_{t\lm\infty} \|B_t(x)\|/a^t =0 \ \forall \ x\in \clm$,
then $\exists \ \hat{a} < a$ and $C$ such that $\|B_t(x)\| < C \hat{a}
^t \ \forall \ x\in \ \clm, \ t\ge 0$.
\newline
(b) Further, suppose that $a \le 1$ and $\lim_{t\lm\infty}\|A_t(x)\|
\|B_t(x)\|^b = 0 \ \forall \ x\in \ \clm$.  Then $\exists \ 
\hat{b}<b$ and $C$ such that $\|A_t(x)\|\|B_t(x)\|^{\hat{b}} < C \ 
\forall \ x \in \ \clm, \ t\ge 0$.
\newline
(c) If $\kk(x) < a \le 1$ and $\theta(x) < b \ \forall \ x\in \ 
\clm$, then $\lim_{t\lm \infty} \|B_t(x)\|=0$ and
$\lim_{t\lm \infty} \|A_t(x)\| \|B_t(x)\|^b=0$ uniformly on $\clm$
\newline
(d) $\kk$ and $\theta$ attain their suprema on $M$.
\end{lemma}
We can now state a corollary to Theorem \ref{thm2}.
\begin{cor}
\label{semi-group}
Assume the generalized Lyapunov numbers of the semi-group $S_t$
satisfy $\kk(p) < 1$ and $\theta(p)<1$ for all $p\in \overline{M'}$,
and that $S_t$ leaves the submanifold $\clm$ negatively invariant.
Further, suppose that there is a $C^1$ semi-group $S^\kappa_t$ which is
$\kappa$-$C^1$ close to $S_t$ on $V$ in some bounded time interval.
Then for $\kappa$ sufficiently
small, there is a submanifold $\overline{M^\kappa}$ which is left
negatively invariant to $S^\kappa_t$, and converges to $\clm$ as
$\kappa \lm 0$ in $C^0$.
\end{cor}
\paragraph{Proof.}
By Lemma \ref{uniformity}, we may choose $T$ sufficiently large so that
with $S=S_T$, Conditions \ref{condB} and \ref{condA}
%and \ref{invert} 
are satisfied.  Thus it follows from Theorem \ref{thm1} that the graph
transform associated with $S_T$ has a fixed point in $\sd$ which we
call $\uk$.

Using the same argument as in Lemma \ref{lemma3}, we see that for
small $t>0$, $h^{-1}(S^\kappa_t(h(Im(\uk)))) \cap 
N^\eps|_{\cup U^{J-3}_i}$ is 
the graph of a section $u^\kappa_t \in \sd$.  Since $h(Im( \uk))$ is
negatively invariant to $S^\kappa_t$, and $S^\kappa_t$ commutes with
$S^\kappa_T$, we see that $u^\kappa_t = \uk$. We can then repeat the
argument for all $t$.
\quad $\blacksquare$

%%%%%%%%%%%%%%%%%%%%%%%%%%%%%%%%%%%%%%%%%%%%%%%%%%%%%%%%%%%%
\section{The PDE}
\label{sectionpde}
\setcounter{equation}{0}

We consider a PDE that may be expressed as an evolution
equation on a separable Hilbert space $\h$ in the form
\begin{equation}
\begin{array}{c}
\frac{du}{dt}+Au+R(u)=0\;,\\
u(0)=u_0\;.
\end{array}
\label{pde}
\end{equation}
We denote the inner product in $\h$ by $(\cdot,\cdot)$ and norm 
$|\cdot|^2= (\cdot,\cdot)$.
We assume that $A$ is a densely defined sectorial linear operator 
with compact inverse.
% eigenvalues $\{\lambda_i\}$ and associated
%eigenfunctions $\{\varphi_i\}$. 
Thus it is possible to choose $\zeta \ge 0$ such that all eigenvalues of
${\tilde A}:=A+\zeta Id$ have strictly positive real part. 
For $\alpha>0$ we define  
${\tilde A}^\alpha=({\tilde A}^{-\alpha})^{-1}$, where
\begin{equation}
{\tilde A}^{-\alpha}=\frac{1}{\Gamma(\alpha)}\int_{0}^{\infty}
t^{\alpha-1}e^{-{\tilde A}t} dt\;.
 \label{001}
\end{equation}
We denote by $\d({\tilde A}^\alpha)$ the domain of ${\tilde A}^\alpha$
(see \cite{He}, \cite{Pazy}).
For $\alpha=0$, we define ${\tilde A}^0=Id$. Then
$\d({\tilde A}^\alpha)$ is a Hilbert space with 
the inner product $({\tilde A}^\alpha u,{\tilde A}^\alpha v)$
and norm $|u|_\alpha=|{\tilde A}^{\alpha}u|$ for 
all $u,v\in \d({\tilde A}^\alpha)$.

The operator $A$ generates an analytic semi-group $L(t)$.
We assume that the nonlinear term $R\in C^1(\d({\tilde A}^{\gamma}),
\d({\tilde A}^{\gamma-\beta}))$ and satisfies
\begin{equation}
\begin{array}{c}
|R(u)|_{\gamma -\beta}\leq C+M(\rho)\hspace{.25in}\forall
u\in \d(A^\gamma):  |u|_\gamma\leq\rho \\
\|DR(u)\|_{op} \equiv \sup_{|v|_\gamma=1} |DR(u)v|_{\gamma-\beta}  \leq 
M(\rho)\hspace{.25in}\forall
u\in \d({\tilde A}^\gamma):  |u|_\gamma\leq\rho\;,
\end{array}
\label{3.2'}
\end{equation}
with $\gamma\geq 0$ and $\beta\in[0,1)$ and where  $C>0$ and
$M(\cdot):\rrr^+\rightarrow \rrr^+$ is a given monotonically-increasing
continuous function with $M(0)=0$. 
Furthermore, for all $u\in\d(\tagm)$ about which we linearize $R$,
we require that $A+DR(u)$ be a closed sectorial operator on $\h$ 
(see \cite{He} for the definitions and properties of sectorial 
operators) with domain
$\d(\ta)$ which equals $\d(A)$ since $\ta$ remains a closed densely
defined operator on $\h$.  In fact, for $\gamma \in [0,1]$ we have 
\begin{lemma}
\label{closed}
Let $A$ be a sectorial operator on $\h$ with $\sigma(A)>0$ with $A^{-1}$
bounded, and let $B$ be a linear operator on $\h$ with $\d(A)\subset
\d(B)$ and which is $A$-bounded
in the sense that for some $K>0$,
$|A^{\gamma-\beta}Bu| < K |A^\gamma u|$ for all $u\in
\d(A)$ and all $\beta \in [0,\gamma)$.  Then $A+B$ is closed in $\h$ 
with domain $\d(A)=\d(A+B)$.
\end{lemma}
\paragraph{Proof.}
Let $x_n\lm x$ and $(A+B)x_n \lm y$ in $\h$.  Then $x_n + A^{-1}Bx_n
\lm A^{-1}y$ so that $\lim_{n\lm\infty}Bx_n = y -Ax$.  Hence, $y-Ax\in
\h$ and $y\in \h$ so that $Ax\in\h$ and thus $x\in \d(A)$.  Furthermore,
$x_n\lm x$ in $\d(A)$ since $-Ax_n = Bx_n -y$ and $A^{-1}(\lim_{n\lm
\infty}Bx_n -y) = x$ in $\h$.

It now suffices to show that $Bx_n \lm Bx$ in $\h$ in order to ensure
that the graph of $A+B$ is closed with domain $\d(A)$.  We simply
estimate
\begin{eqnarray*}
|Bx_n - Bx| &\le & \|A^{\beta-\gamma}\|_{op}|A^{\gamma-\beta}B(x_n-x)|\\
            &\le & K\cdot |A^\gamma(x_n-x)|\\
            &\le & \mbox{using a Poincare-type inequality}\\
            &\le & \tilde{K}\cdot |Ax_n-Ax|,
\end{eqnarray*}
and this completes the proof.
\quad $\blacksquare$

We remark that although it is shown in Henry [\cite{He}, proof of 
Theorem 1.3.2] that $A+B$ has a bounded resolvent on the resolvent set,
it is not shown that $A+B$ is closed.  Therefore, in conjunction with
Lemma \ref{closed}, we have that $A+B$ is sectorial.  In the case
that the spectrum of $A$ does not have strictly positive real part, we
analyze $\ta$ defined as above.

Finally, we assume that the solution operator of \R{pde}, denoted by
$S(t):\d(\tagm)\lm\d(\tagm)$, for $t\ge 0$.  In particular, $S(t)$
is a strongly continuous semi-group, such that for each $t>0$, $S(t)$
is continuously differentiable on $\d(\tagm)$, and its Fr\'{e}chet
derivative is bounded on any bounded set.

Partial differential equations that satisfy the assumptions of this
section include the two-dimensional Navier-Stokes equations, 
the Kuramoto-Sivashinsky equation, the complex Ginzburg-Landau equations, 
certain reaction-diffusion equations, as well as many other PDEs.

%%%%%%%%%%%%%%%%%%%%%%%%%%%%%%%%%%%%%%%%%%%%%%%%%%%%%%%%%%%%%%%%%%%%%%%%

\vspace{.15in}
\noindent {\bf Perturbation of the PDE.}

We consider $C^1$ perturbations of the mapping $G:\d(\tagm)\lm 
\d(\tagm)$ obtained as the time-$\tau$ map of the semi-group of \R{pde}, 
$G(u):=S(\tau)u$ for some fixed $\tau>0$. We denote the perturbed 
mapping by $G_h$ and require that for any $R>0$, there
exist a $K=K(R)\ge0$ such that

\vspace{.15in}
\noindent {\bf Assumptions ${\bf G_h}$}
\begin{itemize}
\item $|G(v)-G_h(v)|_\gamma=|E(v)|_\gamma \leq 
K(R)h\hspace{.25in} \forall v\in \b_R(0)$,
\item $\|DG(v)-DG_h(v)\|_{op}=\|DE(v)\|_{op}\leq K(R)h\,\,\,
\forall v\in \b_R(0)$,
where $DE$ is the Fr\'{e}chet differential of $E$.
\end{itemize}

We remark that $G_h$ may arise from a fully discrete numerical
approximation of the PDE.  For example, let $X^{h,\gamma}$ be a subspace
of $\d(\tagm)$.  In the case of a finite element-type approximation,
$X^{h,\gamma}$  is a finite-dimensional subspace spanned by a family of
polynomials.  Let $P_h$ be the projection taking $\d(\tagm)$ onto $X^{h,
\gamma}$, and let $\tilde{G}_h$ be the finite element approximation
to the map $G$. In order to define our approximation on the same space,
we define $G_h= \tilde{G}_h \circ P_h$.  Thus, $G_h:\d(\tagm)\lm
\d(\tagm)$, and for the largest element diameter $h$ taken sufficiently
small, satisfies Assumptions $G_h$ by standard approximation arguments.

%%%%%%%%%%%%%%%%%%%%%%%%%%%%%%%%%%%%%%%%%%%%%%%%%%%%%%%%%%%%%%%
\section{Unstable Manifolds of hyperbolic fixed points}
\label{section4}
\subsection{Preliminaries}
\setcounter{equation}{0}

Assume that $\bar{u}$ is a hyperbolic fixed point of the PDE \R{pde}.
We first prove the existence of a local $C^1$ unstable manifold of
$\bar{u}$ which by definition, is an overflowing invariant manifold of
\R{pde}.  Then, we show that the global unstable manifold
of the hyperbolic fixed point, formed by evolving the local unstable
manifold forward in time, persists under $C^1$ perturbations of the 
semi-group generated by \R{pde}.

We define the new variable
\[
v(t)=u(t)-\bar{u}\;.
\]
For some fixed $\tau>0$ we define the map  $G(v)=S(\tau)(v+\bu)-\bu$.
Thus, we have translated the hyperbolic fixed point to zero.
Hence, $v(t)$ satisfies 
\begin{equation}
\begin{array}{c}
\frac{dv}{dt}+Cv+F(v)=0 \\
v(0)=v_0=u_0-\bar{u},
\end{array}
\nonumber
\end{equation}
where
\begin{equation}
C=A+D\!R(\bar{u}),
\hspace{.2in}F(v)=R(v+\bar{u})-R(\bar{u})-D\!R(\bar{u})v\;,
\nonumber
\end{equation}
and where $D\!R(\bar{u})v$ is the Fr\'{e}chet derivative of $R$ at
$\bar{u}$ in the direction $v$.

Since $\bu$ is a hyperbolic stationary solution, our assumptions
ensure that 
the operator $C$ is a densely defined sectorial operator with a
spectral gap about the imaginary axis.
We let $\sigma_1(C)$ and $\sigma_2 (C)$ be the spectral sets of $C$ 
associated with the spectrum of $C$ in the left- and right-half planes,
respectively.  Furthermore, we require the cardinality of $\sigma_1$ to
be finite.   We associate the spectral projector $P:\h \lm \h$ with
$\sigma_1(C)$ and $Q:\h \lm \h$ with $\sigma_2$ and so we have that
$\h = P\h \oplus Q\h.$
Using the result of Lemma \ref{closed}, we may use
standard techniques to verify for the Navier-Stokes
equation below, that for any $a\in (0,1)$, there exists $\tau>0$ such
that
\begin{eqnarray}
\begin{array}{c}
\|e^{Ct}v\| \le a \|v\| \quad \forall t \ge \tau,\;\forall v  \in
P\h \cap \d({\tilde A}^{\gamma})
 \\
\|e^{-Ct}v\| \le a \|v\| \quad \forall t \ge \tau,\;\forall v 
\in Q\h \cap \d({\tilde A}^{\gamma})\;,
\end{array}
\label{Cbounds}
\end{eqnarray}
where $\|\cdot\|=|\cdot|_\gamma$ for some $\gamma\in(0,1]$.
We remark that the estimates (\ref{3.2'}) of the nonlinear term dictate
which $\gamma$ to choose for the space $\d({\tilde A}^{\gamma})$. 

Let
\begin{equation}
G(v)=Lv+N(v), \label{Gmap}
\end{equation}
where
\[
Lv=
e^{-C\tau}v, \hspace{.2in} N(v)=\int_{0}^{\tau}e^{-C(\tau-s)}
F(S(s)v)ds\;, \]
with $S(s)$ the semi-group of \R{pde}. 
It is clear that
$D\!N(0)=0$ and that $L$ has no eigenvalues on the unit circle.
Because of the smoothing properties of the operator $\exp(-C\tau)$,
the map $G: \d({\tilde A}^{\gamma})\rightarrow  \d({\tilde
A}^{\gamma})$.

We denote the ball of radius $\rho$ centered about the origin in 
$\d({\tilde A} ^\gamma)$ by
%\[
%\|L\|_{op}:=\sup_{\scriptstyle u\in \d({\tilde A}^{\gamma})\atop\scriptstyle
%|A^\gamma u|=1}|A^\gamma Lu|
%=\sup_{\scriptstyle u\in \d({\tilde A}^{\gamma})\atop\scriptstyle
%\|u\|=1}\|Lu\|\;,
%\]
\[
\b_\rho(0):=\{v\in \d({\tilde A}^{\gamma}) :\,\,\|v\|\leq\rho\}\;.
\]
In general, our requirements on the nonlinear term and on its
linearization are sufficient to ensure that
$N$ $\in$  $C^1(\d({\tilde A}^{\gamma})$,$\d({\tilde A}^{\gamma}))$ 
(see Lemma \ref{needlemma} for details)
and that positive constants $K_1$, $K_2$ exist such that
\begin{equation}
\begin{array}{c}
\|\rr N(v)\|\leq K_1 M(\rho)\rho \hspace{.2in}\forall v\in \b_\rho(0) \\
\|\rr D\!N(v)\|_{op}:=\sup_{\|\eta\|=1}\|\rr D\!N(v)\eta\|
\leq K_2M(\rho) \hspace{.2in}\forall v\in \b_\rho(0)\;, 
\end{array}
\label{Nbounds}
\end{equation}
where $D\!N$ is the Fr\'{e}chet differential of $N$,
$M(\rho)\rightarrow0$ as $\rho\rightarrow0$ and $\rr$ is either
$Id$, $P$, $Q$.

Just as in the proof of  Corollary \ref{semi-group}, we may 
show that the existence of a $C^1$ local unstable manifold
for the map implies that the manifold is
an unstable manifold for the continuous equations. This is quite
advantageous when studying the behavior of overflowing sets of \R{pde}
under numerical approximation.

\subsection{Local Overflowing Manifolds for ${\bf G}$}
To show that there exists a $C^1$ local unstable manifold
of $\bar{u}$ we consider the linear system $Lv=e^{-C\tau}v$
for some fixed $\tau>0$ so that \R{Cbounds} holds. By assumption, the
space $P\h \cap \d({\tilde A}^{\gamma})$ is invariant under $L$ and 
$P\h \cap  \d({\tilde A}^{\gamma})\bigcap \b_\rho(0)$ is a $C^\infty$
overflowing manifold for $L$ for any $\rho>0$.

\begin{lemma} Equation \R{pde} has a $C^1$ overflowing invariant
manifold $W^{\bu,\rho}$ in neighborhood of $\bar{u}$.
\label{gover}
\end{lemma}

\paragraph{Proof.} 
Matching our notation from Section \ref{theory}, we set our 
infinite-dimensional phase space $X = \h$.
We set our unperturbed
mapping $S=L$, the linear map, and our perturbed mapping $\sk=G$.
By construction, the subspace $P\h\cap\d(\tagm)$ is invariant to $L$, 
so we set $M= P\h\cap\d(\tagm) \cap \b_\rho(0)$ for $\rho>0$.  
Then $\clm$ is a
negatively invariant (to $L$) flat manifold of class $C^\infty$, and
we have the splitting $T_x\h \cong\h = T_xM\oplus N_x$ for all $x\in M$, 
where we
may identify $T_x M$ with $P\h\cap\d(\tagm)$, and $N$ with $Q\h\cap
\d (\tagm)$.  Our tubular neighborhood map is simply $h(p,v) = p+v$.
 
We must verify Conditions \ref{condB} and \ref{condA}
with  the constant $C$ in those conditions set to $1$.  
Since $L$ is linear operator independent of $p\in P\h\cap\d(\tagm)$,
for Condition \ref{condB}, we must show that
\begin{equation}
\|Q L v \| \le \vartheta_1 \|v\| \ \ \forall \  v \in
Q\h\cap\d(\tagm).
\nonumber
\end{equation}
This is precisely (\ref{Cbounds}) with $\vartheta_1=a$. For Condition
\ref{condA}, we need $\|P L^{-1}|_{P\h\cap\d{\tagm}} \| \cdot \|QL\| \le
\vartheta_1$, which is trivially satisfied by (\ref{Cbounds}) since
$a^2 < a$.  
%For the geometry under consideration, the Lipschitz sections
%of  the normal bundle $\sd$ are simply $\sd = \{ \Psi \in C^1
%(\clm, Q\h\cap\d(\tagm): \sup_{p\in M} \le \eps, \|\Psi(p_1)-
%\Psi(p_2)\| \le
%\|p_1-p_2\| \ \ \forall \ \ p_1,p_2 \in \clm\}$. 
%Notice that Condition \ref{invert} is trivially satisfied here even
%when the invariant manifold is not compact. This is due to the fact that
%$PL(p,\Psi(p))$ has bounded inverse on $M$ due to
%the commutivity of $P$ and $L$ and  the independence of $L$ on elements
%of $\sd$, together with \R{Cbounds}.

Thus, we have shown that (\ref{Gmap}) has an overflowing invariant
manifold, and appeal to Corollary \ref{semi-group} to yield the result.
\quad $\blacksquare$

\subsection{Persistence of global unstable manifolds}

Herein, we study the behavior of the set 
\begin{equation}
 W^{\bu}:=\cup_{m\in\n} G^m(W^{\bu,\rho})
\label{gun}
\end{equation}
under $C^1$ perturbations of the mapping $G$.
We assume that $W^{\bu}$ is relatively compact and contained in $\b_R(0)$
for some $R>0$.

We define the asymmetric Hausdorff semi-distance  by
\begin{equation}
\mbox{dist}(A,B):=\sup_{a\in A}\inf_{b\in B}\|a-b\|\;,
\label{hausdorff}
\end{equation}
where $A$, $B$ are subsets of $\d({\tilde A}^{\gamma})$. We recall
that dist$(A,B)=0$ if and only if $\bar{A}\subset\bar{B}$.

\begin{thm} 
\label{thm3}
Let $W^{\bu,\rho}$ be the local unstable manifold of
the hyperbolic fixed point $\bar{u}$ of the mapping $G$ given by
Lemma \ref{gover}, and let $G_h$ be a $C^1$ perturbation of $G$ 
satisfying Assumptions $G_h$.  Then, there exists a $C^1$ manifold
$W^{\bu,\rho}_{h}$, overflowing invariant to $G_h$.
Moreover, if $W^{\bu}$ is relatively
compact in $\d({\tilde A}^{\gamma})$, the set
\begin{equation}
W^{\bu}_{h}:= \cup_{m\in\n} G_{h}^m (W^{\bu,\rho}_{h})
\label{ghun}
\end{equation}
satisfies $dist(W^{\bu},W^{\bu}_{h})\rightarrow0$ as $h\rightarrow0$.
\end{thm}

\vspace{.1in}
\noindent{\bf Remark.}
The dynamics on $W^{\bu}$ may be different than on $W^{\bu}_{h}$
and $W^{\bu}$ may be a proper subset of $W^{\bu}_{h}$ in the limit 
$h\rightarrow0$. We illustrate this with an example given in
the appendix.

\vspace{.1in}
\paragraph{Proof.}  
We first show that the mapping $G_h$ has a $C^1$ negatively invariant
manifold.  To do so, we choose $h>0$ sufficiently small  and take
$\rho$ in Lemma \ref{gover} smaller if necessary so that $G_h$ is
in a small enough $C^1$ neighborhood of $L$, so that the conclusions
of Lemma \ref{gover} hold for $G_h$.  From this, we may conclude that
$G_h$ has a unique $C^1$ negatively invariant manifold $\wuph$.

Next, we show that $\wuph \lm \wup$ as $h\lm 0$ in $C^0$.  To do so,
we now consider $G_h$ to be a perturbation of $G$ rather than $L$.
We let $X=\h$ (or $\d(\tagm)$), the unperturbed
mapping $S=G$, and the perturbed mapping $\sk=G_h$.  
We let $\clm
=\wup$ which by Lemma \ref{gover} is negatively invariant with respect
to $G$. From,
Theorem \ref{thm2}, $\wup$ is $C^1$ close to $P\h\cap\d(\tagm)$, hence
the normal and tangent bundles of these two manifolds are $C^0$ close,
respectively.  This, together with the continuity of the norm and the
$C^1$ closeness of $L$ and $G$, ensure the satisfaction of Conditions
\ref{condB} and \ref{condA} with $\vartheta_1 =a'$, where 
$\rho$ is taken smaller if necessary in Lemma \ref{gover} to satisfy
$a'\in [a,1)$.  
% Since $G$ can be made arbitrarily $C^1$ close to $L$,
% and $L$ trivially satisfies Condition \ref{invert}, it is clear that
% $G$ satisfies Condition \ref{invert} as well (see the remark below
% \R{4}).
Thus, applying Theorem \ref{thm1}, $G_h$ has a negatively invariant
manifold $\bar{W}$ which converges to $\wup$ as $h\lm 0$ in $C^0$.

With $\sd$ set to the Lipschitz sections over $\wup$, in order
to show that $\bar{W}= \wuph$, we need only prove that $\wuph \in
\sd$, and then appeal to the  uniqueness of the contraction mapping
theorem.  Let $\nup$ be the normal bundle of $\wup$, and let
$\Puph:P\h\cap\d(\tagm)\cap\b_\rho(0)\lm\d(\tagm)$ be the $C^1$ section 
identified with the overflowing invariant manifold $\wuph$ by the tubular
neighborhood map.  Then, each fiber of $\nup$ is a submanifold of 
$\d(\tagm)$.  If there does not exist an element of $\sd$ which
can be identified with $\wuph$, then  there is a $p\in \wup$ such 
that $\Puph$ is transversal to $\nup|_p$, contradicting the fact that
$\wuph$ is $C^1$ close to $\wup$.  Thus, $\wuph$ is a section over
$\wup$. The fact that this section is Lipschitz follows from  the
$C^1$ closeness of $\wup$ and $P\h\cap\d(\tagm)\cap \b_\rho(0)$.

Next, we prove the lower semi-continuity of $\wup$.  Since $\wu$ is
relatively compact in $\d(\tagm)$, for any $\epsilon'>0$, there exists
$u^i \in \wu$, $i=1,...,N$, for some $N < \infty$ such that
$\wu \subset\cup_{i=1}^N \b_{u^i} (\epsilon'/2)$.  Hence, for any $u\in
\wu$, $u_h\in\d(\tagm)$,
\begin{equation}
\|u-u_h\|\leq\frac{\epsilon'}{2}+\|u^i-u_h\|,
\label{temp1}
\end{equation}
for some $i$, $1\leq i\leq N$. 

To proceed with the argument we need to find
$u^{i}_{h}\in \wuh$ close to $u^i$ in $\d(\tagm)$.  Choose $n$ large
enough so that $G^{-n}u_i\in W^{\bu,\rho}$ for all $1\leq i\leq N$. 
Due to the global triviality of the unperturbed manifold $\wup$, we 
may once again express the sections of the normal bundle of this 
manifold as $\s = \{ \Psi \in C^1 (\clm, Q\h\cap\d(\tagm) 
\cap \Gamma(N^\eps) \}$.   (The argument that we just gave
using the underlying $C^\infty$ flat manifold shows that $\s$ is
not empty.) Let the image
of $\Phi_h$ in $\s$ be the fixed point of the mapping $G_h$ so that
$\wuph = {\text h} \circ \Phi_h|_{\overline{\wup}}$, and define
$u_h^i = G^n_h(G^{-n}(u^i), \Phi_h(G^{-n}(u^i)))$, where h is the
tubular neighborhood map.

We obtain from Assumptions $G_h$ that for any $w,z \in \b_R(0)$
\begin{equation}
\|G(w)-G_h(z)\|\leq K_L\|w-z\|+K(R)h\;,
\label{lipball}
\end{equation}
where we may assume without loss of generality that $K_L>1$.
If $w\in \wup$, and $z\in\d(\tagm)$ is defined by
$(Pz,Qz)=(w,\Phi_h(w))$, then $\|w-z\| \lm 0$ as 
$h\lm 0$ since $\wup$
and $\wuph$ can be made arbitrarily $C^0$ close.
By assumption $\|G^n(w)\|<R$ for all $n\geq0$, $w\in\wup$.
With this choice of $w$ and  $z$, \R{lipball} 
implies that $\|G_h(z)\|<R$, by shrinking $h$ if necessary.
Thus we may continue to apply Assumptions $G_h$.
By iterating this procedure we find that
\begin{eqnarray}
\|G^n(w)-G^{n}_{h}(z)\|
&\leq& K^{n}_{L}\|w-z\|+Kh\sum_{i=0}^{n-1}K^{i}_{L} \nonumber \\
&\leq& K^{n}_{L}\|w-z\|+\left(\frac{K_{L}^{n}}{K_L-1}\right) Kh\;.
\label{4.9}
\end{eqnarray}
We apply this last inequality with $w=G^{-n}u^i$ and $z\in\d(\tagm)$
defined by
$(Pz,Qz)=(w,\Phi_h(w))$.  As before $\|w-z\| \lm 0$ as 
$h\lm 0$. Hence, for $h$ sufficiently small, using (\ref{4.9}) we
obtain that $ \|u^i-u_{h}^{i}\| \le \epsilon'/2$.

Therefore, \R{temp1} yields
\begin{equation}
\inf_{u^h\in W^{\bar{u}}_h}\|u-u^h\|\leq 
\frac{\epsilon'}{2}+\|u^i-u_{h}^{i}\|\leq \epsilon'\;.
\end{equation}
We conclude that dist$(W^{\bu},W^h)\leq\epsilon'$
and so dist$(W^{\bu},W^h)\rightarrow0$ as $h\rightarrow0$.  
\quad $\blacksquare$

The previous theorem only makes explicit use of the properties
of the overflowing invariant manifold, but does not use the properties
of hyperbolic fixed points {\em per se}; consequently, this result
generalizes to any overflowing invariant manifold. In particular, we may
combine it with a result of Hale, Lin \& Raugel \cite{HLR} to obtain

\begin{thm} 
Let $\aaa$ be the global attractor for \R{pde} which we assume to be
compact in $\d(\tagm)$.
Further, suppose that the union of overflowing invariant manifolds which
satisfy Conditions \ref{condB}, \ref{condA},
%\ref{invert} 
and their global unstable manifolds, formed by evolving the
overflowing manifold forward in time, is dense in $\aaa$.
Then, every $C^1$ perturbation of
\R{pde}, in the sense of Assumptions $G_h$, has an attractor $\aaa^h$
that satisfies $\max[dist(\aaa,\aaa^h),
dist(\aaa^h,\aaa)]\rightarrow0$ as $h\rightarrow0$.
\label{hausthm}
\end{thm}

\paragraph{Proof.} 
That $dist(\aaa,\aaa^h)\rightarrow0$ as
$h\rightarrow0$ follows the proof of the previous theorem.
That $dist(\aaa^h,\aaa)\rightarrow0$ as $h\rightarrow0$ is given in
\cite{HLR}.
\quad $\blacksquare$

\vspace{.1in}
An example of such a system is a gradient system (see \cite{Hale}).

%%%%%%%%%%%%%%%%%%%%%%%%%%%%%%%%%%%%%%%%%%%%%%%%%%%%%%%%%%%%%%%%
\section{Application to the Navier-Stokes equations.}
In this section we verify that the Navier-Stokes  equations (NSE)
and their fully discrete approximations satisfy the assumptions of 
Sections \ref{sectionpde} and \ref{section4}. 
We will consider the two-dimensional NSE for a viscous incompressible
fluid in a bounded open simply-connected domain $\Omega$ with boundary
$\partial \Omega$ of class $C^2$.
We consider no slip boundary conditions.
The governing equations are
\begin{equation}
\left.\begin{array}{c}
\frac{\partial u}{\partial t} +(u\cdot\nabla)u=-\nabla p 
+\nu \Delta u +f
\,\,\,\,\mbox{in}\,\,\,\Omega\times(0,\infty) \\
\nabla\cdot u=0\,\,\,\,\mbox{in}\,\,\,\Omega\times(0,\infty)\\
u|_{\partial\Omega}=0,
\end{array}\right\} \label{NSE}
\end{equation}
where $f=f(x)$, the volume force, and $\nu>0$, the kinematic viscosity,
are given. We denote by 
$u=u(x,t)$ the velocity vector, and $p=p(x,t)$ the pressure which are
the unknowns.

\subsection{Estimates on the nonlinear term and exponential dichotomies}
Following the notation in \cite{CF},
we set
\begin{eqnarray*}
\vv&=&\{v: v\in (C^{\infty}_{0})^2, \,\,\,\mbox{div}\, v=0\},
\end{eqnarray*}
and define
\[
\begin{array}{c}
\h=\mbox{the closure of}\,\,\, \vv\,\,\, \mbox{in}\,\,\, (L^2(\Omega))^2,\\
\v=\mbox{the closure of}\,\,\, \vv\,\,\, \mbox{in}\,\,\, (H^1(\Omega))^2,
\end{array}
\]
where $L^2(\Omega)=H^0(\Omega)$ and $H^l(\Omega)$, $l=1,2,\ldots$
denotes the usual $L^2$-Sobolev spaces. $\h$ is a Hilbert space
with the $L^2$ inner product and norm.

Let $P_\h$ denote the orthogonal projection in $L^2(\Omega)\times 
L^2(\Omega)$ onto $\h$. We denote by $A$ the Stokes operator
\[
Au=-P_\h\Delta u,
\]
and the bilinear operator
\[
B(u,v)=P_\h((u\cdot\nabla)v)
\]
for all $u,v$ in  $\d(A)=\v\cap (H^2(\Omega)\times H^2(\Omega))$. 
The operator $A$ is a self-adjoint  positive-definite 
operator with  compact inverse. 

We recall some bounds on the bilinear term $B(u,v)$. The first is
\begin{equation}
|B(u,v)|\leq c|u|^{1/2}\,|A^{1/2}u|^{1/2}\,\,|A^{1/2}v|^{1/2}|Av|^{1/2}
\hspace{.3in} \forall u,v,\in \d(A), \label{b2}
\end{equation}
where $c$ will denote an adequate positive constant.
Alternatively, we may use Agmon's inequality
\[
\|u\|_{L^{\infty}(\Omega)}\leq K |u|^{1/2}\,|Au|^{1/2}
\hspace{.3in}\forall u\in \d(A),
\]
to obtain
\begin{equation}
|B(u,v)|\leq K|u|^{1/2}\,|Au|^{1/2}\,\,|A^{1/2}v|
\hspace{.7in}\forall u\in \d(A),
\,\,\,\,v\in V. \label{b6}
\end{equation}

\begin{lemma} Assume that $f\in \v$ and let $S(t)$ be the semi-group
operator of the NSE. Then  $S(t)$ is strongly continuous and for each
$t>0$, $S(t)u$ is continuously differentiable on $\d(A)$ with locally
bounded Fr\'{e}chet derivative.
Furthermore, the nonlinear term of the NSE satisfies (\ref{3.2'}) with 
$\gamma=1$ and $\beta=1/2$, and the global attractor is compact in 
$\d(A)$.
\label{nsenonlinear}
\end{lemma}

\paragraph{Proof.} 
The Cattabriga-Solonnikov-Vorovich-Yudovich theorem (see
\cite{Catta}, \cite{Solon}, \cite{VY}),
yields constants  $C_1, C_2$ depending on $\Omega$
such that 
$C_1\|u\|_{(H^1(\Omega))^2}\leq 
|A^{1/2}u|\leq C_2\|u\|_{(H^1(\Omega))^2}$ for all $u\in \d(A^{1/2})$,
giving us the equivalence of norms.
In addition, the projection 
$P_\h:(H_{0}^{1}(\Omega))^2\lm (H^{1}(\Omega))^2$ is bounded
(see for example \cite{CF}). Hence, we find that for all $u,v \in
{\mathcal V}$
\[
|A^{1/2}B(u,v)|^2\leq K\int_{\Omega}\sum_{i,j,k=1,2}
\left|\frac{\partial u_i}{\partial x_k}
\frac{\partial v_j}{\partial x_i}+u_i\frac{\partial^2v_j}
{\partial x_k\partial x_j}\right|^2d\Omega\;,
\]
and so for all $u,v \in \d(A)$ by a density argument.
We use inequalities of the type \R{b2} on the terms
$\frac{\partial u_i}{\partial x_k}\frac{\partial v_j}{\partial x_i}$
in the sum and  \R{b6} on the terms $u_i\frac{\partial^2v_j}
{\partial x_k\partial x_j}$ to obtain
\begin{eqnarray*}
|A^{1/2}B(u,v)|&\leq& K|A^{1/2}u|^{1/2}\,|Au|^{1/2}\,
|A^{1/2}v|^{1/2}\,|Av|^{1/2}\\
& & + \,\,K|u|^{1/2}\,|Au|^{1/2}\,|Av|^{1/2}
\hspace{.3in} \forall u,v \in \d(A)\;.
\end{eqnarray*}
Then, since $\lambda_1|u|\leq|Au|$,
$\lambda_{1}^{1/2}|A^{1/2}u|\leq |Au|$, where $\lambda_1>0$ is the
first eigenvalue of $A$, we obtain
\begin{equation}
|A^{1/2}B(u,v)|\leq K|Au|\,\,|Av| \hspace{.2in} \forall u,v \in \d(A)\;,
\label{hardnon}
\end{equation}
and so the inequalities of (\ref{3.2'}) are satisfied with 
$R(u)=B(u,u)-f$  and $\gamma=1, \beta=1/2$. 
Further, for any $u_0$ such that $|Au_0|\leq R$, we
have that $|Au(t)|\leq K(R)$ for all $t\geq 0$ and some
$K(R)>0$ (see \cite{CF}, \cite{Lion}, and \cite{T}).

The proof that $S(t)u$ is continuously  differentiable for all $t>0$ is
similar to Lemma 14.3 in \cite{CF}, merely needing modification of the
topology from $\h$ to $\d(A)$, and also shows that
$\|DS(t)|_{\b_{R}(0)}\| \le M(R) < \infty$ for all $R>0$.
Because we have assumed that the force
$f\in \d(A^{1/2})$, it follows that there is an absorbing ball
in $\d(A^{3/2})$, and hence the global attractor is in fact compact 
in $\d(A)$ (\cite{CF}).
\quad $\blacksquare$

Let us suppose that $\bu$ is a stationary solution such that
the operator
\[
Cv=Av+B(\bu,v)+B(v,\bu)
\]
%has a finite number of eigenvalues with negative real part and
has no eigenvalues with zero real part. 
As in the previous section we construct the map $G$ as follows.
We set $v(t)=u(t)-\bu$ and
\[
F(v):=B(v+\bu,v+\bu)-B(\bu,\bu)-DB(\bu,\bu)v=B(v,v).
\]
Then 
\begin{eqnarray}
G(v)&=&Lv+N(v) \label{NSEG}\\
&=&S(\tau)(v+\bu)-\bu,\nonumber
\end{eqnarray}
where 
\vspace{-.15in}
\[
Lv=e^{-C\tau},\hspace{.2in} N(v)=\int_{0}^{\tau}e^{-C(\tau-s)}F(v(s))ds.
\]
\begin{lemma} 
The operator $C$ is sectorial and $G:\d(A)\lm \d(A)$.
Moreover, there exists a $\tau>0$ such that \R{Cbounds}, \R{Nbounds} are
verified.
\label{needlemma}
\end{lemma}
\paragraph{Proof.} From Lemma \ref{closed}, we obtain that $C$ is closed in $\h$ with 
domain $D(C)=D(A)$.  In order to show that $C$ is sectorial we need
prove that $(C-A)A^{-\alpha}$
is a bounded operator on $\h$ for some $\alpha\in[0,1)$ according
to [\cite{He}, Corollary 1.4.5].
Using \R{b2} and \R{b6}, we obtain that
\begin{eqnarray*}
|(C-A)A^{-1/2}v| &\leq& |B(\bu,A^{-1/2}v)|+|B(A^{-1/2}v,\bu)|\\
&\leq& c|A\bu|\,|v|+c|v|\,|A\bu|.
\end{eqnarray*}
That $C$ has compact resolvent and hence  has a discrete spectrum is
given to us by analyticity arguments as in [\cite{CFT}, Lemma 2.1].
Thus, we may appeal to [\cite{He},
Theorem 1.4.8] to obtain the equivalence of the norms $|C_1^\beta 
\cdot|$ and $|A^\beta \cdot|$ for all $\beta \in [0,1]$, where  we
define $C_1=C+cId$ for $c>0$ such that $Re(\sigma(C_1))$.

As in Section 4.1, we associate the spectral projector $P$  with the 
spectral set $\sigma_1=Re(\sigma(C))<0$, and  the projector $Q$ with 
the its complement $\sigma_2=\sigma/\sigma_1$. Here, $\sigma(C)$ is the
spectrum of $C$ and we have the splitting $\h = P\h \oplus Q\h$.

Using, [\cite{He}, Theorems 1.5.3 and 1.5.4], we have that
\begin{eqnarray}
|Ae^{-Ct}v|\leq C_3e^{\alpha t}|Av|\hspace{.2in}\forall v\in \d(A)
,\,\,\,\, t\geq0, \label{5.1}\\
|Ae^{Ct}p|\leq C_3e^{-\gamma t}|Ap|\hspace{.2in}\forall p\in P\h\cap
\d(A) ,\,\,\,\, t\geq0,
\label{5.1'}\\
|Ae^{-Ct}q|\leq C_3t^{-1/2}e^{-\gamma t}|A^{1/2} q|
\hspace{.2in}\forall q\in Q\h\cap \d(A),\,\,\,\,t>0,
\label{henrysave}
\end{eqnarray}
where $-\alpha< Re(\sigma_1(C))$, 
$0<\gamma<Re(\sigma_2)$ and for some $C_3>0$.
We may also write (\ref{henrysave}) as $|Ae^{-Ct}q| \le C_3 t^{-1/2}
e^{-\gamma t} K |Aq|$ for all $q\in Q\h\cap\d(A)$ and $t>0$, by 
appealing to the positivity of $A$.
Thus given $a<1$ there is a $\tau>0$ such that 
$\max\{C_3e^{-\gamma \tau},C_3\tau^{-1/2}e^{-\gamma \tau}K\}\leq a$
and \R{Cbounds} are verified. 

We show that
\begin{equation}
|Ae^{-Ct}v|\leq K t^{-1/2}|A^{1/2}v|\hspace{.2in}\forall v\in \d(A)
\ \ \ 0<t<\tau.
\label{henry7}
\end{equation}
Indeed, for $v\in \d(A)$ and with $p=Pv$ and $q=Qv$ we have that
\begin{eqnarray*}
|C_1 e^{-Ct}v| & \le & |C_1 e^{-Ct}p| + |C_1 e^{-Ct}q| \\
 & \le & \mbox{(using \R{5.1'} and \R{henrysave} with the equivalence 
of norms)}\\
 & \le & C_3e^{-\gamma t} K |C_1 p| + C_3t^{-1/2}e^{-\gamma t} K 
|C_1^{1/2} q| \\ 
 & \le & \mbox{(since $P\h\cap\d(A)$ is finite-dimensional)} \\
 & \le & C_3e^{-\gamma t} K |C_1^{1/2} p| + C_3t^{-1/2}e^{-\gamma t} K 
|C_1^{1/2} q| \\ 
 & \le & \mbox{(and for $t \in (0,\tau]$ )}\\
 & \le & C_3t^{-1/2}e^{-\gamma t} K |C_1^{1/2} p| + 
C_3t^{-1/2}e^{-\gamma t} K |C_1^{1/2} q| \\ 
 & \le & \mbox{(using the commutivity of $C$ with the projectors
$P$ and $Q$)}\\
 & \le & C_3t^{-1/2}e^{-\gamma t} K |C_1^{1/2} v|.
\end{eqnarray*}
($K$ is a generic positive constant.)
Hence, \R{henry7} follows,  and we find that
\begin{equation}
|AN(v)|\leq \int_{0}^{\tau}\frac{K}{(\tau-s)^{1/2}}|A^{1/2}F(v(s))|ds
\leq K\tau^{1/2}\max_{0\leq t\leq \tau}|Av(t)|^2.
\nonumber
\end{equation}
Standard energy estimates yield an estimate of the form
$|Av(t)|\leq K(\tau)|Av(0)|$ for $0\leq t\leq \tau$
({\em c.f.} the analysis of [\cite{T}, Theorem 10.2] and 
\cite{CF}) and the first inequality in \R{Nbounds} 
is satisfied since the projectors $P$ and $Q$ are bounded on $\h$.

We also have that
\[
DF(v)w=DB(v+\bu)w-DB(\bu)w=B(v,w)+B(w,v).
\]
Thus in the same way one finds that
\[
\|ADN(v)\|_{op}=\sup_{|A\eta|=1}|ADN(v)\eta|\leq K|Av|\;,
\]
and the second inequality in \R{Nbounds} is satisfied. \quad $\blacksquare$

\subsection{Perturbations of the NSE}
To illustrate the verification of Assumptions $G_h$ we consider
the fully discrete approximation of the NSE given by
\[
u^{h}_{n+1}=(Id+\dt\nu A)^{-1}[u^{h}_{n}+\dt 
P_N(f- B(u^{h}_{n},u^{h}_{n}))]\;,
\]
where $u^{h}_{N}\in P_N\h$ and where $P_N$ is the projection
onto the first $N$ eigenfunctions of $A$.
With $T$ chosen so that \R{Cbounds} holds, we define the map
\begin{eqnarray}
u^{h}_{n}&=&
(Id+\dt\nu A)^{-n}u^{h}_{0}+\dt\sum_{j=1}^{n}(Id+\dt\nu A)^{j-n-1}
P_N[f-B(u^{h}_{j-1},u^{h}_{j-1})]\nonumber\\
&=&\tilde{G}_{h}(u^{h}_{0}),
\label{timedis}
\end{eqnarray}
where $n\dt=T$.

\begin{lemma} The map defined by 
\[
G_h(v)=\tilde{G}_{h}(P_N(v+\bu))-\bu
\]
approximates the map $G$, given by \R{NSEG}, in the $C^1$ norm
on any bounded set.  In particular, Assumptions $G_h$ are satisfied.
\end{lemma}

\noindent {\it Sketch of Proof}. It is sufficient to prove that the maps
$G(u)=S(T)u$ and $G_h(u)=\tilde{G}^{h}(P_Nu)$ are close in the $C^1$
norm on $\b_R(0)$ for any $R>0$ and for $h$ sufficiently small
({\em i.e.} $N$ sufficient large and $\dt$ sufficiently small).
As shown in \cite{JS1}, $G$ and $G_h$ are 
continuously close. The proof is straightforward but tedious.

To show $C^1$ convergence is not difficult;
one copies the proof for the continuous convergence,
only applied to the linearization of the equations (as above).
That is, one shows that 
\[
DG_h(u^{h}_{0})\mu^{h}_{0}=
(Id+\Delta t\nu A)^{-n}P_N\mu^{h}_{0}-\sum_{j=1}^{n}
(Id+\Delta t\nu A)^{j-n-1}P_ND\!B(P_Nu_{j-1}^{h})\mu_{j-1}^{h}
\]
approximates $DG(u)$ in precisely the same way 
one shows that \R{timedis} and $S(T)$ are close.
One finds that for all $u\in \b_R(0)$ (we omit the details)
\[
|A(G(u)-G_h(u))|+
\|A(DG(u)-DG_h(u))\|_{op}\leq K(R,T)\dt
+\frac{K(R,T)}{\lambda_{N+1}^{1/2}}\;.
\]

As shown in Lemma \ref{nsenonlinear}
the global attractor is in fact compact in $\d(A)$.
Since $W^{\bu}$ is part of the global attractor, it is
relatively compact in $\d(A)$.
Thus we may restate Theorem~\ref{hausthm} for the map $G_h(u)$, \R{timedis},
with $h=\max\{\dt,\lambda_{N+1}^{-1/2}\}$:

\begin{thm} Any overflowing invariant manifold satisfying the
conditions of Section 2 is captured by the fully discrete
approximation, \R{timedis}, $h$ sufficiently small.
Moreover, the global unstable manifold of this overflowing manifold,
formed by evolving the overflowing manifold forward in time is
lower semi-continuous in the Hausdorff semi-distance sense at $h=0$.
\end{thm}

%%%%%%%%%%%%%%%%%%%%%%%%%%%%%%%%%%%%%%%%%%%%%%%%%%%%%%%%%%%%%%%
\section{Inertial Manifolds}
\setcounter{equation}{0}

In this section we study the persistence of inertial manifolds
associated with a nonlinear evolution equation 
in a Hilbert space $\h$, taking the form
\begin{equation}
\begin{array}{c}
\frac{du}{dt}+Au+R(u)=0\ \\
u(0)=u_0\;,
\end{array}
\label{pde1}
\end{equation}
under $C^1$ perturbation of the semi-group.  The assumptions  on
$A$ and $R(u)$ are the same as for \R{pde}, but in addition, we assume
that \R{pde1} is dissipative and that $A$ is self-adjoint, meaning that 
$u(t)\in \b_{R_0}(0)$ for all $t \ge T^*(u_0)$,
where $R_0>0$ is independent of $u_0$.

We assume that the eigenvalues of $A$, $0<\lambda_{1}
\leq \lambda_{2}\ldots\leq \lambda_{j}\rightarrow\infty$, 
repeated with their multiplicities,
satisfy for any positive $K_3$, $K_4$
\begin{equation}
\lambda_m\geq K_3,\hspace{.35in}\lambda_{m+1}-\lambda_{m}
\geq K_4 \lambda_{m+1}^{\beta}
\label{gap}
\end{equation}
for some $m\geq m_0$, $m_0 > 0$. We define $P_m$ to be 
the projection onto the 
span of the generalized eigenfunctions of the operator $A$ 
corresponding to the first $m$ eigenvalues, and set $Q_m$ to be the
projector corresponding to the complementary spectral set.
PDEs that are known to satisfy the above assumptions include
the Kuramoto-Sivashinsky equation, the complex Ginzburg-Landau 
equations, certain reaction-diffusion equations, as well as other PDEs.
(See the references that are listed in the introduction.)
The above assumptions are sufficient to guarantee that \R{pde1} has an
inertial manifold. Indeed, we recall

\begin{thm} Suppose that \R{gap} is satisfied. Then for any
$R\geq R_0$, $\epsilon>0$, $\delta\leq 1$ there exists a
sufficiently large $m$
so that \R{pde1} has an inertial manifold representable as the graph
of a $C^1$ function $\Phi:P_m\d(A^\gamma)\lm Q_m\d(A^\gamma)$ 
contained in 
$\b_R(0)$. Moreover, $\Phi$ has the properties 
\[
\sup_{ \|p\|\leq R}\|\Phi(p)\|\leq\epsilon\;, \hspace{.3in}
\sup_{\|p\|\leq R}\|D\Phi(p)\|_{op}\leq\delta\;.
\]
\label{imthm}
\end{thm}

Again, $S(t)$ is the semi-group operator of (\ref{pde1}) and we set 
$G(u)=S(\tau)u$ for some $\tau>0$.  As usual, we decompose $G(u)=
Lu+N(u)$, where
\begin{equation}
Lu:=e^{-A\tau}u, 
\quad N(u):=-\int_{0}^{\tau}L(\tau-s)R(S(s)u)ds\;. \label{plus}
\end{equation}
Clearly, the inertial manifold of \R{pde1} is the inertial manifold for
$G$ as well, so to study its persistence, we consider $C^1$ 
perturbations $G_h$ of the map $G$, satisfying Assumptions $G_h$ of 
Section \ref{sectionpde}.

In general, the $C^1$ perturbed mapping $G_h$ are not dissipative
even though the map $G$ is assumed to be so.
Thus although Theorem \ref{thm2} states that $G_h$ possesses an
inflowing invariant manifold, it may not be an inertial manifold.
For this reason, we will restrict our attention to the truncated 
map ${\tilde G_{h}}$ which we now define.

Let $\theta: \rrr^+\mapsto[0,1]$ be a fixed $C^1$ function such
that $\theta(x)=1$ for $0\leq x\leq2$, $\theta(x)=0$ for
$x\geq 4$, and $|\theta'(x)|\leq 2$ for $x\geq 0$.
Define $\theta_R(x)=\theta(x/R^2)$.
We consider the map
\begin{eqnarray}
u&=&{\tilde G_{h}}(v):=Lv+N(v)-\theta_{R_{0}}(\|v\|^2)E(v)\nonumber\\
&:=& Lv+N^h(v)\;,
\label{nhdef}
\end{eqnarray}
where $L$ is given (\ref{plus}).
The map ${\tilde G_{h}}$ agrees with the map $G_h$ inside
the ball $\b_{\sqrt{2}R_0}(0)$ and it is dissipative.  Moreover, if
we require Assumptions $G_h$ to hold on $\b_{2R_0}(0)$, then
\begin{equation}
\|G(v)-{\tilde G_{h}}(v)\|\leq K(R)h,\hspace{.25in}
\|DG(v)-D{\tilde G_{h}}(v)\|_{op}\leq K(R)h\,\,\,\,\,
\forall v\in \d(A^\gamma)\;.
\label{globalc1}
\end{equation}

Unlike the negatively invariant manifolds that we examined in previous 
sections, an inertial manifold is an example of an
inflowing or positively invariant manifold.  For this reason, we must 
redefine our
definitions of the linear operators $A$ and $B$ given in Conditions
\ref{condB} and \ref{condA}.  We set
\begin{equation}
A(p)\equiv  DG|_{M'}(p):T_pM' \lm T_{G(p)}M' \text{  and  }
B(p)\equiv  Q \circ DG(p)|_{N}:N_p \lm N_{G(p)}.
\nonumber
\end{equation}

\begin{thm} Suppose that the map ${\tilde G_{h}}$ satisfies
\R{globalc1} and the the map $G$ has an inertial manifold
as given in Theorem~\ref{imthm}. Then for $h$ sufficiently 
small, the map ${\tilde G_{h}}$ has an inertial manifold  
in $\b_R(0)$ which converges to the true manifold as $h$ tends to zero.
\end{thm}

\paragraph{Proof.} We set $X=\d(A^\gamma)$, our unperturbed mapping
$S=G|_{\b_{R_0}(0)}$, and our perturbed mapping $\sk={\tilde G_{h}}$.
The unperturbed inflowing manifold $M$ is the inertial manifold for
$G$. In order to show that ${\tilde G_{h}}$ has an inflowing invariant
manifold we appeal to Theorem \ref{thm1}, and again verify
Conditions \ref{condB} and \ref{condA}.
%, and \ref{invert}. 
Lemmas 5.4 and 5.5 of \cite{RT} shows that
that the inertial manifold $M$ is normally hyperbolic.
Keeping the notation in Conditions \ref{condB} and \ref{condA}, these
lemmas yield
\[
\|A(p)\| \leq Ke^{(\lambda_m+K_2\lambda_{m}^{\gamma-\beta})\tau}\;,
\]
\[
\|B(p)\| \leq
Ke^{-(\lambda_{m+1}-2K_2(1+\delta)\lambda_{m+1}^{\gamma-\beta})\tau}
\]
for all $p\in P_m\d(A^\gamma)\cap\b_{R_0}(0)$.
Hence, \R{gap} with $\tau=\lambda_{m+1}^{-1}$ verifies
Conditions \ref{condB} and  \ref{condA}.

%Finally, we verify Condition \ref{invert}. Since the unperturbed
%inertial manifold is globally trivial, use a single chart
%$U_1^{J-3}\equiv P_m\d(A^\gamma)\cap\b_{R_0}(0), \sigma_1)$, and
%since $X$ is a
%vector space, we let the tubular neighborhood map $h(p,q)=p+q$.
%We may then let $(G^1,G^2)$ be the local coordinate representation
%of the map $G$. In this case $(G^1,G^2)=(P_mG,Q_mG)$.
%We need to show that $D_1G^1(x,y)$ has a bounded inverse for
%in a some $X$-neighborhood of $P_m\d(A^\gamma)\cap\b_{R_0}(0)$.
%We have from \R{plus} that
%\begin{eqnarray*}
%D_1G^1(x,y)&=&e^{-A\tau}+\int_{0}^{\tau}P_mD_1R(S(s)(x,y)ds\\
%&=&e^{-A\tau}\Big(Id+e^{A\tau}\int_{0}^{\tau}P_mD_1R(S(s)(x,y)
%ds\Big)\;.
%\end{eqnarray*}
%Using the bounds on $R(u)$ in Section \ref{sectionpde} and the fact that
%$\| A^\beta e^{-tA}\|_{\l(\h,\h)} \le Kt^{-\beta}$ (see \cite{He} for
%example), we obtain for $u\in \b_R(0)$
%\begin{eqnarray*}
%\left\|e^{A\tau}\int_{0}^{\tau}P_m D_1R(S(s)(x,y)\right\|_{op}\
%&\leq& 
%\left\|e^{A\tau}\int_{0}^{\tau}P_m DR(S(s)(x,y)\right\|_{op}\ \\
%&\leq& 
%e^{\lambda_m\tau}\int_{0}^{\tau}\frac{RK_2M(\rho)}
%{(\tau-s)^\beta}ds\\
%&\leq& e^{\lambda_n\tau}\frac{\tau^{1-\beta}}{1-\beta}RK_2M(\rho)
%\;.
%\end{eqnarray*}
%If we choose $\tau=\lambda_{m+1}^{-1}$, we make this last
%inequality less than one (by increasing $m$ if necessary) and so
%$D_1G^1$ has a bounded inverse. 

Thus, Theorem \ref{thm1} applies and the mapping ${\tilde G_{h}}$ 
has an inflowing invariant manifold that converges to $M$ in $C^0$.
This inflowing invariant manifold is also an inertial manifold for the
mapping ${\tilde G_{h}}$.
\quad $\blacksquare$

In \cite{JS1} and \cite{JSTi}, the perturbed manifold is constructed
off of the linear space $P_m\d(A^\gamma)$. As a consequence,
a restriction on the size of the time step $\tau$ in the map of
equation (\ref{plus})
is required. In addition, the dimension of the inertial manifold 
may be required to increase to insure that certain estimates hold.
We have avoided both of these issues by constructing the perturbed 
manifold off of the (inflowing invariant) inertial manifold.

\vspace{.1in}

If the map $G_h$ is dissipative with absorbing set $\b_{R_0}(0)$, then
we may take $R=2R_0$ in Theorem~\ref{imthm} and merely check
Assumptions $G_h$ for $R=2R_0$. Furthermore, the
inertial manifold for ${\tilde G_{h}}$ is an inertial manifold for the
map $G_h$ inside $\b_{R_0}(0)$. A more generic case is given by

\begin{lemma} Suppose that the ball $\b_{R_0}(0)$ attracts orbits of
\R{pde1} exponentially with uniform rate. Then for any $R_h>R_0$ there 
exists an 
$R_h>R_{h}^{c}>R_0$ such that for $h$ sufficiently small the ball
$\b_{R_h}(0)$ is positively invariant under $G_h$ and 
$\b_{R^{c}_{h}}(0)$ attracts orbits starting in $\b_{R_h}(0)$ of $G_h$
exponentially with uniform rate. Moreover, for such $h$ the inertial 
manifold for
${\tilde G_{h}}$ is an inertial manifold for the map $G_h$ restricted to
$\b_{R_h}(0)$.
\end{lemma}

\paragraph{Proof.}
By assumption we may suppose that $\|u(t)\|\leq 
R_0-(R-R_0)e^{-\lambda t}$ for
all $u(t)$ solving \R{pde1} with $\|u_0\|\leq R$,
and for some $\lambda>0$. Assume that Assumptions $G_h$ hold for the
given $R_h>R_0$. 
Let $h$ be taken sufficiently small so that
$K(R_h)h+R_0+(R_h-R_0)e^{-\lambda\tau}<R_h$ where $K(R_h)$ given in
Assumptions $G_h$. 
Then we have that for $u\in \b_R(0)$ with  $R\leq R_h$ 
\[
\|G_h(u)\|\leq \|G_h(u)-G(u)+G(u)\|
\leq K(R_h)h+R_0+(R-R_0)e^{-\lambda\tau}\;.
\]
There is a minimum value, $R=R^c_h$, for which the above inequality is
less than $R$.  It is given by
\[
R_{h}^{c}=\frac{R_0(1-e^{-\lambda\tau})+K(R_h)h}{1-e^{-\lambda\tau}}\;.
\]
Of course $R_{h}^{c}\rightarrow R_0$ as $h\rightarrow 0$. We choose $h$
sufficiently small so that $R_{h}^{c}\leq\sqrt{2}R_0$ ( the map
${\tilde G_{h}}$ and $G_h$ agree inside $\b_{\sqrt{2}R_0}(0)$). Then for
any $u\in \b_{R_h}(0)$ the inequality above shows that $G_h(u)\in 
\b_{R_h}(0)$
and $G_{h}^{n}(u)$ converges exponentially to $\b_{R_{h}^{c}}(0)$ as
$n\rightarrow \infty$.
\quad $\blacksquare$

\vspace{.1in}
Persistence of inertial manifolds under numerical approximation has
been studied elsewhere. In \cite{FST} and  \cite{FSTi}, the persistence
of an inertial manifold under a Galerkin approximation of the 
underlying PDE is studied, while \cite{DG2} studies the behavior
of the inertial manifold under a time discretization.
In \cite{FTi} and \cite{Jones1}, finite differencing approximations
are studied, and in \cite{JS1} more general maps are
considered, as well as finite element approximations and their time 
discretizations. 

%%%%%%%%%%%%%%%%%%%%%%%%%%%%%%%%%%%%%%%%%%%%%%%%%%%%%%%%%%%%%%%
\section{Appendix}
\setcounter{equation}{0}
We provide an example of the remark we made following Theorem 
\ref{thm3}, in which the sequence of unstable manifolds $\wuh$ do not 
converge to the true dynamics of $\wu$ in the limit as $h$ tends to 
zero.

We consider the vector field expressed in polar coordinates
\[
\dot{r}=(r-1)(r\sin^2\frac{\theta}{2}-1)+
h r\sin^2\frac{\theta}{2}
\]
\[
\dot{\theta}=\sin(\theta) \; ,
\]
where $h\geq0$. Notice that for $h=0$ there is a stationary solution
at $r=1$, $\theta=\pi$ which is not hyperbolic. For $h>0$ this stationary
solution disappears. Notice also that the vector field with $h>0$
is a $C^1$ perturbation of the vector field with $h=0$.

\begin{figure}[ht]
\epsfxsize=0.9\textwidth
\centerline{\epsfbox{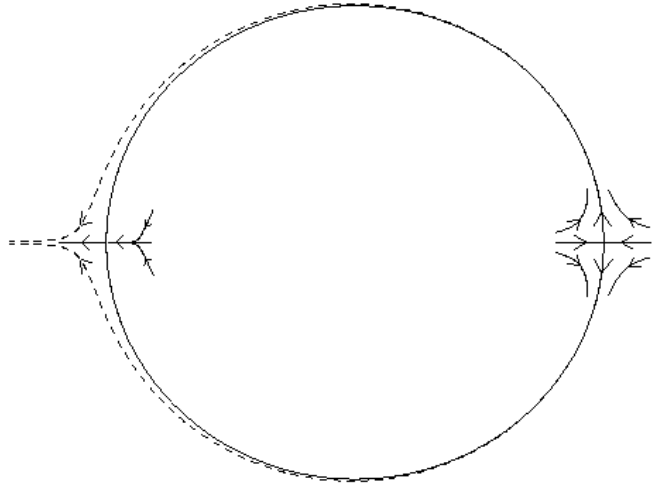}}
\caption{\label{figure}
{\it\small Lower Semi-continuity of approximate unstable manifold to the
true unstable manifold .} \normalsize}
\end{figure}

The dashed line in Figure \ref{figure} 
represents the flow of the perturbed system ($h>0$); the solid lines 
represent the 
unperturbed system ($h=0$). The unstable manifold of
the hyperbolic point at $r=1$, $\theta=0$ of the perturbed system
now approximates the unstable manifold of the non hyperbolic 
stationary solution of the unperturbed system as well as
the unstable manifold of the unperturbed hyperbolic fixed point. By 
preparing the above vector field these manifolds remain bounded.

Notice that the unperturbed unstable manifold
of the hyperbolic point $r=1$, $\theta=0$ is
a proper subset of the perturbed unstable manifold in the
limit as $h\rightarrow0$.

\noindent{\bf Acknowledgements.}
The authors are grateful to Edriss Titi for numerous discussions
and helpful ideas.  The authors  would also like to thank Peter
Constantin for suggestions on estimating the nonlinear term of
the NSE, as well Bernard Minster and Len Margolin for their 
encouragement and interest.
S.S. gratefully acknowledges the support of the Cecil H. and Ida M.
Green Foundation and the Institute of Geophysics and Planetary
Physics at Scripps Institution of Oceanography.
D.A.J. gratefully acknowledges the support of the Institute for
Geophysics and Planetary Physics (IGPP) and the Center for
Nonlinear Studies (CNLS) at Los Alamos National Laboratory.
This work  was partially supported by the Department of Energy 
``Computer Hardware, Advanced Mathematics, Model Physics" (CHAMMP) 
research program as part of the U.S. Global Change Research Program.

\end{document}